\documentclass[11pt, reqno]{article}

\usepackage[utf8]{inputenc}
\usepackage[T1]{fontenc}

\usepackage{setspace}
\usepackage[margin=1in]{geometry}
\usepackage{parskip}
\usepackage[sc]{mathpazo}
\usepackage[T1]{eulervm}
\usepackage[scaled]{helvet}

\usepackage{amsmath,amsthm, amssymb, amsfonts}
\usepackage{mathtools}
\usepackage{graphicx}
\usepackage[hidelinks]{hyperref}

\newcommand{\arxiv}[1]{\href{https://arxiv.org/abs/#1}{\texttt{ArXiv:#1}}}

\theoremstyle{plain}
\newtheorem{thm}{Theorem}
\newtheorem{cor}{Corollary}[section]
\newtheorem{lem}{Lemma}[section]
\newtheorem{prop}{Proposition}[section]

\theoremstyle{definition}

\newtheorem{exm}{Example}[section]

\numberwithin{equation}{section}

\newcommand{\E}[1]{\mathbf{E}\left[#1\right]}

\newcommand{\ind}[1]{\mathbf{1}_{\{ #1 \}}}
\newcommand{\dt}[1]{\mathrm{det}\left(#1\right)}
\newcommand{\intz}[1]{\frac{1}{2 \pi \mathbf{i}}\oint \limits_{#1}}

\newcommand{\Z}{\mathbb{Z}}
\newcommand{\R}{\mathbb{R}}
\newcommand{\C}{\mathbb{C}}
\newcommand{\W}{\mathbb{W}}

\newcommand{\cD}{\mathcal{D}}

\title{On the duality between particles and polymers}
\author{Mustazee Rahman\thanks{\textsc{Department of Mathematical Sciences, Durham University}. \textit{Email}: \texttt{mustazee@gmail.com}}}
\date{}

\begin{document}

\maketitle

\begin{abstract}
\noindent We explore the connection between tasep-like interacting particle systems and last passage percolation type polymer models,
focusing on three models: Geometric, Exponential and Brownian last passage percolation and their associated tasep particle systems. We explain how formulas for certain natural observables in last passage percolation translate to formulas for tasep, by going through a notion of "duality". In turn, we obtain determinant formulas for last passage percolation with a deterministic boundary and for tasep with a deterministic first particle trajectory.
\end{abstract}

\section{Introduction} \label{sec:intro}
\subsection{Duality}
There is a close connection between certain interacting particle systems with one-sided interaction (such as tasep) and certain directed 2-d polymer models (such as last passage percolation). This connection is made through a variational formula which, for tasep-like systems, expresses their evolution in terms of the initial condition and a metric semigroup law of the associated polymer model. It is the probabilistic analogue of the Hopf-Lax-Oleinik formula, which gives the solution of certain Hamilton-Jacobi equations (such as Burgers' equation) as a variational principle involving the initial condition and a certain Legendre dual. As such, a good word to describe this connection would be "duality". Duality is like a dictionary between two languages. It need not be a bijection.

It is well known that the hydrodynamics of tasep (which is the abbreviation of ``totally asymmetric simple exclusion process") is given by Burgers' equation \cite{Liggett}. Perhaps the earliest realization of this duality within probability is due to Rost \cite{Ro}, who saw that the hydrodynamics of tasep could be used to establish limit shapes for Exponential last passage percolation. Hydrodynamics and limit shapes are laws of large numbers. What is remarkable about this duality between particles and polymers is that it extends all the way to fluctuations. This is the connection between the KPZ fixed point \cite{MQR} and the directed landscape \cite{DOV} as explained in \cite{DV, NQR, RV}. Duality provides a coupling between these two objects in terms of a variational principle which translates many properties between them.

Let's consider this duality from the viewpoint of Brownian last passage percolation and Brownian tasep. Brownian last passage percolation is defined in terms of a family $W_1,W_2,\ldots$ of independent, standard Brownian motions. It is a stochastic process defined over the domain $\{(s,n;t,m) \in (\R \times \Z)^2: 0\leq s \leq t; 1 \leq n \leq m\}$.
The last passage value from $(s,n) \to (t,m)$ is
\begin{equation} \label{eqn:blppintro}
    L(s,n;t,m) = \max_{s = t_{n-1} \leq t_n \leq t_{n+1} \leq \cdots \leq t_m = t} \sum_{k=n}^m W_k(t_k) - W_k(t_{k-1}).
\end{equation}
Brownian tasep (also called reflected Brownian motions \cite{WFS}) is an interacting particle system
\begin{equation} \label{eqn:Btasepintro} X(t) = (X_1(t) \geq X_2(t) \geq X_3(t) \geq \cdots)\end{equation}
where $t \geq 0$ and $X_n(t) \in \R$. Its evolution is described through a family $W^*_1, W^*_2, \ldots$ of independent Brownian motions. Given an initial condition $X_1(0) \geq X_2(0) \geq \cdots$, $W^*_n$ is a Brownian motion started from $X_n(0)$. Particle 1 has trajectory $X_1 = W^*_1$. Particle 2 moves according to $W^*_2$ reflected to the left off $X_1$. In general, particle $n$ has trajectory $X_n$ which is $W^*_n$ reflected to the left off $X_{n-1}$. The reflection is in the sense of Skorokhod (see below).
The relation between Brownian last passage percolation and Brownian tasep is
\begin{equation} \label{eqn:bbduality}
    X_n(t) = \min_{1 \leq k \leq n} \{X_k(0) - L(0,k;t,n)\}.
\end{equation}
The equality here is not just in law for two stochastic processes. It is more of less a deterministic identity. One can simply \emph{define} the left-side in terms of the right-side; the Brownian motions $W_k$ are transmuted deterministically to $W^*_k$. Duality produces a dictionary to express probabilities for $X$ in terms of $L$ and vice-versa. In doing so, it ``permutes" the roles of space, time and noise. For such a thing to work, space, time and noise should occupy the same dimension, ergo, (1+1)-dimensional models. In this article, we try to understand duality between tasep-like systems and last passage percolation type models. Along the way, we find theorems that come about from the language alone; so we spend some time developing the language itself.

\subsection{Skorokhod reflection}
Let $f_1, f_2 : [0,\infty) \to \R$ be continuous functions with $f_2(0) \geq f_1(0)$. The Skorokhod reflection, upwards, of $f_2$ off $f_1$ is the function
\begin{equation} \label{eqn:reflectiondef}
    f_2 \odot f_1 (t) = \max \{f_2(t), \max_{s \in [0,t]} f_1(s) + f_2(t)-f_2(s)\}
\end{equation}
for $t \geq 0$. Visually, the graph of $f_2 \odot f_1$ is the graph of $f_2$ pushed upwards from the graph of $f_1$. When $f_2(0) = f_1(0) = 0$,
the maximum in \eqref{eqn:reflectiondef} is obtained by the second argument: $(f_2 \odot f_1)(t) = \max_{s \in [0,t]} f_1(s) + f_2(t)-f_2(s)$.

Now let $B(t)$ be a standard Brownian motion and consider $B \odot b$, where $b(t)$ is a continuous function with $b(0) = 0$.
We have the following determinant formula for the finite dimensional laws of $B \odot b$ as the $n=1$ case of Theorem \ref{thm:blppintro} below.

\begin{prop} \label{prop:Btimesb}
    Let $B(t)$ be a standard Brownian motion and $b: [0,\infty) \to \R$ be continuous with $b(0)=0$.
    Let $0 < t_1 < t_2 < \ldots < t_k$. Then,
    $$\Pr( B \odot b(t_i) \leq a_i; 1 \leq i \leq k) = \dt{I + Q + R}_{L^2(\{1,\ldots,k\} \times \R)}.$$
    The determinant is the Fredholm determinant understood through its series expansion. The integral kernels $Q$ and $R$ are as follows.
    \begin{align*}
        Q(i,u;j,v) &= \frac{1}{\sqrt{2 \pi (t_j-t_i)}}e^{- \frac{(u-v)^2}{2 (t_j-t_i)}} \, \ind{t_i < t_j; u \geq a_i; v \geq a_j}; \\
        R(i,u;j,v) &= \int_{-\infty}^{\infty} dz\, \partial_v \Pr(\tau \leq t_j, W(t_j) \in dv \mid W(0)=z) \, \ind{u \geq a_i; v \geq a_j}.
    \end{align*}
    Here $W : [0,\infty) \to \R$ is a standard Brownian motion and $\tau = \inf \{t \geq 0: W(t) \leq b(t)\}$ is a hitting time.
\end{prop}
Notice that $R$ has rank one and $Q$ is strictly upper triangular in terms of its $k \times k$ block structure. Therefore,
\begin{align*}
\dt{I + Q + R} &= \dt{I+Q}( 1 + \mathrm{tr}((I+Q)^{-1}R)) \\
&= 1 + \mathrm{tr}((I+Q)^{-1}R) \\
&= 1 + \sum_{i=0}^{k-1} (-1)^i \mathrm{tr}(Q^i R).
\end{align*}
In the last line we have used the fact that $Q^k = 0$ to write $(I+Q)^{-1} = 1 - Q + Q^2 \cdots + (-1)^{k-1}Q^{k-1}$.

In the classical setting of Skorokhod reflection, one takes $f_1 \equiv 0$. If $X(t) = f \odot 0(t)$, then
$X(t) = f(t) + L(t)$ where $L(t) = \max_{0 \leq s \leq t} f(s)_{-}$. The process $L$ is like a local time of $X$ at zero;
it is non-decreasing and may increase only at times when $X(t) = 0$ \cite{RY}. If $f$ is a standard Brownian motion $B$
then $X \stackrel{d}{=} |B|$. Skorokhod reflection arises naturally in contexts involving reflected or coalescing Brownian motions;
see, for instance, \cite{STW} for a discussion and also invariance properties of Brownian motion under reflection.
Proposition \ref{prop:Btimesb} appears to be the first instance of a formula for the reflection of Brownian motion off a deterministic function.

\subsection{Brownian last passage percolation and Brownian tasep}
Let $b:[0,\infty) \to \R$ be a continuous function and let $W_1, W_2, \ldots$ be a family of independent, standard Brownian motions.
Brownian last passage percolation with a boundary $b$ is the following process $G$ defined on the set $\{(t,n)\in \R \times \Z: t\geq 0, n \geq 0\}$.
\begin{equation} \label{eqn:blppwithb}
    G(b; t,n) = \max_{0 \leq t_0 \leq t_1 \leq \cdots \leq t_n = t} b(t_0) + \sum_{k=1}^n W_k(t_k) - W_k(t_{k-1}) = \max_{0 \leq t_0 \leq t} b(t_0) + L(t_0,1;t,n)
\end{equation}
where $L$ is from \eqref{eqn:blppintro}.
Visually, $G(b;t,0) = b(t)$, $G(b;t,1)$ is the reflection of $W_1$ off $b$, $G(b;t,2)$ is the reflection of $W_2$ off $G(b;t,1)$, and so on.

Classical Brownian last passage percolation considers the so-called narrow wedge initial condition whereby $b(0)=0$ and $b(t) = -\infty$ otherwise.
In this case, $G(b; t,n) = L(0,1;t,n)$. It is by now a well known fact that $L(0,1;1,n)$ has the same law as the largest eigenvalue of an $n\times n$ GUE random matrix \cite{GTW}.
More generally, the process $n \mapsto L(0,1;1,n)$ has the law of the largest eigenvalue of the minors of an infinite GUE random matrix \cite{Bar}.

\subsubsection{A formula for the law of Brownian last passage percolation with boundary}
We determine the finite dimensional laws of $t \mapsto G(b;t,n)$.
This is stated in Theorem \ref{thm:Blpp} of Section \ref{sec:brown}, and we restate it here.

\begin{thm} \label{thm:blppintro}
Let $n \geq 1$. Let $0 < t_1 < t_2 < \ldots <  t_k$ and $a_1, \ldots, a_k$ be real numbers.
$$Pr(G(b; t_i,n) \leq a_i; 1 \leq i \leq k) = \det(I - \chi_a K_B \chi_a)_{L^2(\{1,\ldots, k\} \times \mathbb{R})}$$
with $\chi_a(i,z) = \ind{z \geq a_i}$ and the integral kernel $K_B$ is as follows.
\begin{equation}
    K_B(i, u; j, v) = - e^{\frac{(t_j-t_i)}{2}\partial^2}(u,v) \, \mathbf{1}_{t_i < t_j} + S_{n,t_i} \cdot S^{hypo(b)}_{n,t_j}(u,v).
\end{equation}
The kernels
\begin{align*}
    e^{\frac{t}{2}\partial^2}(u,v) &= \frac{1}{\sqrt{2\pi t}} e^{-\frac{(u-v)^2}{2t}} \quad t > 0,\\
    S_{n,t}(u,v) &= \frac{t^{(n-1)/2}}{(n-1)!} H_{n-1}\big(\frac{u-v}{\sqrt{t}}\big) \quad t > 0,\\
    S^{hypo(b)}_{n,t}(u,v) & = (-\partial_v)^n \Pr(\tau \leq t, W(t) \in dv \mid W(0)=u) \quad t > 0.
\end{align*}
Here $W$ is a standard Brownian motion, $H_n$ is the $n$-th Hermite polynomial and $\tau = \inf \{t \geq 0: W(t) \leq b(t)\}$.
\end{thm}

\subsubsection{The narrow wedge boundary}
Consider the boundary $b(0)=0$ and $b(t)=-\infty$ otherwise. This is not continuous but the formula in Theorem \ref{thm:blppintro} holds by approximating $b$ with the continuous functions $t \mapsto -\ell t$ and taking $\ell \to \infty$. Indeed, the hitting times $\tau_{\ell}$ associated to the boundary $-\ell t$ converges almost surely and monotonically to the hitting time of $b$, from which one can deduce that the associated kernels $K_B$ converge in trace norm. Furthermore, in \eqref{eqn:blppwithb}, the process $G(-\ell s; t,n)$ converges to $G(b;t,n)$ which can be deduced from the fact that a Brownian motion $B(s)$ satisfies $|B(s)| = O(s^{\frac{1}{2}+})$ in probability.

For the narrow wedge boundary $b$, the hitting time $\tau = \infty$ if $W(0) > 0$ and $\tau = 0$ if $W(0) \leq 0$.
As a result,
$$S^{hypo(b)}_{n,t}(u,v) = (-\partial_v)^n e^{\frac{t}{2}\partial^2}(u,v) \mathbf{1}_{u \leq 0} = t^{-n/2} e^{\frac{t}{2} \partial^2}(u,v) \, H_n \big ( \frac{u-v}{\sqrt{t}} \big) \mathbf{1}_{u \leq 0}.$$

The process $t \mapsto G(b;t,n) = L(0,1;t,n)$ has an interpretation in terms of non-colliding Brownian motions.
Suppose $(\hat{B}_1 < \hat{B}_2 < \cdots < \hat{B}_n)$ are $n$ independent Brownian motions conditioned not to intersect in the sense of Doob's $h$-transform with harmonic function $h(x_1, \ldots, x_n) = \prod_{i < j}\, (x_j-x_i)$ on the domain $\W = \{(x_1,x_2,\ldots, x_n) \in \R^n: x_1 < x_2 < \ldots < x_n\}$. Then $$\hat{B}_n(t) \stackrel{d}{=} L(0,1;t,n)$$ as a process in $t$ \cite{OY}. Furthermore, $\hat{B}_n$ has the law of the trajectory of the top particle among $n$ particles performing GUE Dyson Brownian motion \cite{Gr}. The laws of these processes, as given by Theorem \ref{thm:blppintro}, is fairly well understood now; see \cite{JoDPG}.

\subsubsection{The flat boundary and random matrices}
Consider the boundary $b \equiv 0$. Then
$$G(0;t,n) = \max_{s \in [0,t]}\; L(s,1;t,n).$$
This process is studied in \cite{BFPSW} (it is the process $Y_n$ there).
The authors prove distributional identities for $Y_n(t) = G(0;t,n)$ and $Z_n(t) = L(0,1;t,n)$
in terms of Dyson Brownian motion with reflecting boundaries. These, in turn, are connected to random matrices as we explain.

An $n$-particle Dyson Brownian motion of type $C$ is a stochastic process $X^C(t) = (X^C_1(t), \ldots, X^C_n(t))$ taking values in a Weyl chamber of type $C$:
$\W_C = \{(x_1, \ldots, x_n) \in \R^n: 0 < x_1 < x_2 \cdots < x_n \}$. It evolves according to the stochastic differential equations
$$ dX^C_i(t) = dB_i(t) + \frac{1}{X^C_t(t)} dt + \sum_{j: i \neq j} \left ( \frac{1}{X^C_i(t) - X^C_j(t)} + \frac{1}{X^C_i(t) + X^C_j(t)} \right ) dt.$$
This process models the eigenvalues of a random matrix ensemble \cite{Gr, KT}. Consider the Lie algebra $\mathfrak{so}_{2n+1}$ consisting of $(2n+1) \times (2n+1)$ skew-symmetric matrices.
The eigenvalues of a matrix in $\mathfrak{so}_{2n+1}$ take the form $0, \pm \mathbf{i} \lambda_1, \ldots, \pm \mathbf{i} \lambda_{n}$ with $0 \leq \lambda_1 \leq \cdots \leq \lambda_{n}$.
One can define a Brownian motion over $\mathfrak{so}_{2n+1}$ for which the eigenvalues $\lambda_i(t)$ evolve according to $X_i^C(t)$ started from $X^C_i(0) = 0$.
Much like GUE Dyson Brownian motion (which is Dyson Brownian motion of type A corresponding to the Lie algebra $\mathfrak{su}_n$), the process $X^C$ can be realized as the Doob $h$-transform of $n$ independent Brownian motions killed when it hits the boundary of $\W_C$ and with harmonic function $h_C(x) = \prod_{i=1}^n x_i \cdot \prod_{1 \le i < j \leq n} (x_j^2-x_i^2)$. This can be interpreted as $n$ Brownian motions conditioned to not collide with each other as well as with a wall at zero.

An $n$-particle Dyson Brownian motion of type $D$ is a stochastic process $X^D(t) = (X^D_1(t), \ldots, X^D_n(t))$ taking values in a Weyl chamber of type $D$:
$\W_D = \{(x_1, \ldots, x_m) \in \R^n: |x_1| < x_2 \cdots < x_n \}$. It evolves according to the stochastic differential equations
$$ dX^D_i(t) = dB_i(t) + \sum_{j: i \neq j} \left ( \frac{1}{X^D_i(t) - X^D_j(t)} + \frac{1}{X^D_i(t) + X^D_j(t)} \right ) dt.$$
Now consider the Lie algebra $\mathfrak{so}_{2n}$.The eigenvalues of a matrix in $\mathfrak{so}_{2n}$ take the form $\pm \mathbf{i} \lambda_1, \ldots, \pm \mathbf{i} \lambda_{n}$ with $|\lambda_1| \leq \cdots \leq \lambda_{n}$. The Brownian motion over $\mathfrak{so}_{2n}$ has eigenvalues $\lambda_i(t)$ that evolve according to $X_i^D(t)$ started from zero \cite{Gr, KT}.
This, too, can be obtained as the Doob $h$-transform of $n$ independent Brownian motions killed when it hits the boundary of $\W_D$ and with harmonic function $h_D(x) = \prod_{1 \le i < j \leq n} (x_j^2-x_i^2)$. This can be interpreted as $n$ Brownian motions conditioned to not collide with each other and being reflected off a wall at zero.

In \cite{BFPSW} the Dyson Brownian motion of types $C$ and $D$ are connected to the process $G(0;t,n)$. For type $D$, the authors replace the first component $X^D_1$ by its absolute value,
for which the governing equations become
$$ dX_i^{(D)} = dB_i + \frac{1}{2} \mathbf{1}_{i=1} dL(t) + \sum_{j: j \neq i} \left ( \frac{1}{X^{D}_i - X^{(D)}_j} + \frac{1}{X^{D}_i + X^{(D)}_j} \right) dt.$$
Here $L(t)$ is the local time of $X^{(D)}_1$ at zero. They prove that for every $n \geq 1$,
\begin{align*}
        X^{(C)}_n(t) &\stackrel{d}{=} G(b \equiv 0; t,2n) \\
        X^{(D)}_n(t) &\stackrel{d}{=} G(b \equiv 0; t,2n-1)
\end{align*}

Theorem \ref{thm:blppintro} then produces a determinant formula for the law of a tagged particle in Dyson Brownian motion of types C and D. In Example \ref{exm:flat}
we calculate the corresponding kernel.

One implication of this connection involves the Airy-one process (see \cite{BFPS} for the definition of the Airy-one process). Although we do not pursue it here,
under the so called KPZ scaling the process $x \mapsto n^{1/6}(G(0; (1+ 2xn^{-1/3},n)) -2n^{1/2}-2xn^{1/6})$ should converge to the Airy-one process.
There was a conjecture (see Conjecture 2 in \cite{BFPS}) that the Airy-one process should emerge as the scaling limit of GOE Dyson Brownian motion; this turned out to be false \cite{BFP}.
Nevertheless, we find that the Airy-one process should arise as the KPZ scaling limit of $X_n^{(C/D)}(t)$, which are indeed dynamics on eigenvalues.

The authors of \cite{BFPSW} also observe that for any fixed time $t > 0$,
$$ \max_{s \in [0,t]} \; L(0,1;s,n) \stackrel{d}{=} G(0;t,n).$$
This is an extension (when $n=1$) of the classical reflection principle for Brownian motion, which states that the running maximum of a Brownian motion at a fixed time has the law of reflected Brownian motion. Note that $s \mapsto L(0,1;s,n)$ has the law of the top particle in an $n$-particle GUE Dyson Brownian motion. Thus, the identity shows that the running maximum of GUE Dyson Brownian motion, at fixed times, has the law of Dyson Brownian motion of types C or D. This latter fact is an analogue of (and should imply it in the KPZ scaling limit) the observation of Johansson \cite{JoDPG} that the maximum of the Airy-two process has the GOE Tracy-Widom law.

\subsubsection{Brownian tasep}
Consider Brownian last passage percolation with boundary $b$ and denote by $L_{b}(s,n;t,m)$ the last passage value from $(s,n) \to (t,m)$ according to \eqref{eqn:blppintro} with the driving functions $W_0 = b$ and $W_1, W_2, \ldots$ being independent Brownian motions.
Consider the associated Brownian tasep model as determined by \eqref{eqn:bbduality}.
Particles are now labelled with $0,1,2,\ldots$ and have trajectories $X_0,X_1,X_2,\ldots$. Take the initial positions of the particles to be zero: $X_n(0)=0$ for every $n$. Observe that
$$ X_n(t) = - \max_{0 \leq k \leq n}\, L_b(0,k;t,n) = -L_b(0,0;t,n) = -G(b;t,n).$$
We find that $X_0$ has a deterministic trajectory $-b$, $X_1$ is the reflection (to the left) of $-W_1$ off $X_0$, $X_2$ is the reflection of $-W_2$ off $X_1$, etc.
Theorem \ref{thm:blppintro} then has the following corollary (stated again as Corollary \ref{cor:Blpp} in Section \ref{sec:btasep}).

\begin{cor} \label{cor:btasepintro}
    Consider Brownian tasep started from $X_n(0) = 0$ for $n \geq 0$.
    Let $X_0$ move deterministically according to a continuous function $-b$.
    The other particles move randomly according to independent Brownian motions being reflected to the left off the particle to its right (as noted above).
    Then for $n \geq 1$, times $0 < t_1 <  t_2 <  \ldots < t_k$, and $a_1, \ldots, a_k \in \R$,
    \begin{equation} \label{eqn:btasepformula}
    \Pr(X_n(t_i) \geq a_i; 1 \leq i \leq k) = \dt{I - \chi_{-a} K_B \chi_{-a}}_{L^2(\{1,2,\ldots,k\}\times \R)}
    \end{equation}
    where $\chi_{-a}(i,z) = \ind{z \geq -a_i}$ and the kernel $K_B$ is from Theorem \ref{thm:blppintro}.
\end{cor}

Boundary conditions for last passage percolation become, under duality, the trajectory of the first particle in the corresponding tasep model.
The observables of interest then become the temporal law of a fixed particle trajectory.
For tasep-like models it is also natural to have initial conditions at time zero as in the $X_k(0)$ in \eqref{eqn:bbduality}.
The observables of interest them become the probabilities
$$ \Pr(X_{n_1}(t) \geq a_1, \ldots, X_{n_k}(t) \geq a_k)$$
for some tagged particles $n_1 < n_2 <\cdots < n_k$ at a fixed time $t$. We can think of this as a spatial law at a fixed time.

In \cite{NQR} the authors derive this spatial law of Brownian tasep with arbitrary initial conditions in terms of a Fredholm determinant.
Their argument is based on taking a low density limit of continuous time tasep, and using formulas for the latter as derived in \cite{MQR}.
Brownian tasep is Markovian in time, and Warren \cite{Warren} found a formula for the transition probability of an $N$-particle Brownian tasep.
In Section \ref{sec:warrenformula} we derive Warren's formula from a formula of Johansson for Geometric last passage percolation.
As an application of duality, we explain how one can derive the formula in \cite{NQR} in terms of determinant formulas for Exponential last passage percolation (which are themselves derived from Geometric last passage percolation).

\subsection{Continuous time tasep}
Duality connects continuous time tasep to Exponential last passage percolation. Continuous time tasep is a particle system
$$ X(t) = (X_1(t) > X_2(t) > X_3(t) > \cdots) \quad t \in [0,\infty),$$
where $X_n(t) \in \Z$ is the trajectory of particle $n \geq 1$. Starting from an initial condition $X(0)$ at time zero, each particle attempts to jump one unit rightward at rate 1. A jump is successful if there is no particle blocking the site (the exclusion rule). See \cite{Liggett} for a probabilistic construction of the system.

Suppose we introduce a particle 0 whose trajectory is deterministic in the following way. Choose times $0=\tau_0 < \tau_1 < \tau_2 < \tau_3 < \cdots$ with $\tau_k \to +\infty$ and $X_0(0) > X_1(0)$. Let the particle 0 jump one unit rightward at the times $\tau_k$, that is,
$$ X_0(t) - X_0(0) = \max \{k \geq 0: \tau_k \leq t\}.$$
The other particles move rightward at rate 1 subject to the exclusion rule as before.

Let $X^*(t) = (X^*_0(t) >  X^*_1(t) > \cdots)$ denote the corresponding tasep model with a deterministic first particle trajectory and initial condition $X^*_n(0) = -n$ for $n \geq 0$ (this is called the step initial condition). As an application of duality, we have the following formula for the temporal law of a tagged particle, which is also stated as Corollary \eqref{cor:ctasep} in Section \ref{sec:ctasep}.

\begin{thm} \label{thm:ctasepintro}
    For $n \geq 1$, times $t_1, t_2, \ldots, t_k > 0$, and distinct integers $a_1, \ldots, a_k \geq -n$,
    \begin{equation} \label{eqn:ctasepformula}
    \Pr(X^*_n(t_i) > a_i; 1 \leq i \leq k) = \dt{I - \chi J \chi}_{L^2(\{1,2,\ldots,k\}\times \R)}
    \end{equation}
    where $\chi(i,z) = \ind{z \geq t_i}$ and the kernel $J$ is as follows.
    
    Define the following integral kernels on $\R$ for integers $m,n \geq 1$ and $a \geq -n$.
\begin{align} \label{eqn:kernels}
    Q^m(u,v) &= \frac{(v-u)^{m-1}}{(m-1)!} e^{u-v} \ind{v \geq u}; \\
    S_{n,a}(u,v) &= - \frac{1}{2 \pi \mathbf{i}} \oint_{|w-1|= 1} dw\, e^{(w-1)(v-u)}w^{a+n+1} (1-w)^{-n}\, \ind{v \leq u}; \\
    \bar{S}_{n,a}(u,v) &= \frac{1}{2 \pi \mathbf{i}} \oint_{|w|= 1} dw\, e^{(w-1)(v-u)} w^{-a-n-1} (1-w)^n.
\end{align}
Let $B(m)$ be a random walk with $Exp(1)$ step distribution, whose $m$-step transition probability is given by $Q^m$ above. Define,
$$ \tau = \inf\, \{m = 0,1,2,3, \ldots: B(m) \leq \tau_{m+1}\},$$
and
\begin{equation} \label{eqn:Jhypo}
    S^{hypo}_{n,a}(u,v) = \E{\bar{S}_{n,a-\tau}(B(\tau),v) \ind{\tau \leq a+n} \mid B(0) = u}.
\end{equation}
The integral kernel
\begin{equation} \label{eqn:Jkernel}
    J(i,\cdot;j,\cdot) = - Q^{a_j-a_i} \ind{a_i < a_j} + S_{n,a_i} \cdot S^{hypo}_{n,a_j}.
\end{equation}
\end{thm}

This form of tasep has been recently studied in \cite{BBF} and called tasep with a moving wall. See also \cite{FG, G} for further works on this theme. The authors are interested in the fluctuations of a tagged particle at a single, large time and show that these are governed by the one-point laws of the KPZ fixed point. Theorem \ref{thm:ctasepintro} provides a formula for the temporal law of a tagged particle.

The models that we consider, such as Brownian last passage percolation with a boundary or tasep with a moving wall, are expected to converge to the KPZ fixed point under a joint scaling of space, time and noise (the KPZ scaling; see \cite{MQR}). We do not pursue this scaling limit. One may also take this limit via the variational representation of these models, following \cite{DV}.

\subsection{Further context}
The motivation to study tasep-like particle systems and last passage percolation type polymer models comes from two early sources.
The work of Kardar, Parisi and Zhang \cite{KPZ} introduced the eponymous KPZ equation and predicted the scaling relations under which many particle systems are expected to converge to the recently constructed KPZ fixed point by Matetski, Quastel and Remenik \cite{MQR}. The latter was obtained by finding formulas for the laws of particle trajectories in continuous time tasep, building on the works \cite{Sch, Sas, BFPS}. The formulas we obtain are also based on this general approach.

For polymer models, a problem of Ulam \cite{Ulam} asked about the asymptotic behaviour of the longest increasing subsequence in a large, random permutation. Hammersley \cite{Ham} introduced a key process to solve this problem. Works by Vershik and Kerov \cite{KV}, and by Logan and Shepp \cite{LS}, established the law of large numbers for Ulam's problem via connections to random partitions (see also Aldous and Diaconis \cite{AD} and Sepp\"al\"ainen \cite{Sep} for proofs inspired by particle systems). Baik, Deift and Johansson \cite{BDJ} (see also Borodin, Okounkov and Olshanski \cite{BOO}) proved the central limit theorem and made the connection to random matrices. Johansson \cite{JoSh} made the connection to the last passage percolation models studied here. Dauvergne, Ortmann and Vir\'{a}g \cite{DOV} recently introduced the directed landscape, which is expected to be the scaling limit of many polymer models. The duality between the KPZ fixed point and the directed landscape is shown in \cite{DV, NQR, RV}.

For surveys of these topics and more, see, for instance, \cite{BoGo, CoKPZ, Josurvey, QuKPZ, Zyg} and references therein.

\section*{Acknowledgements}
I am thankful to B\'{a}lint Vir\'{a}g and Kurt Johansson. B\'{a}lint came to visit Durham, and we talked about the use of formulas in probability. Kurt had wondered to me if it's possible to understand formulas appearing in random growth models and exclusion processes from a common viewpoint. These conversations led me to think about the subject of this article. I don't know if it clarifies much of anything, but I did learn a few things. Sometimes, the more you think you understand, the more you realize that you don't know.

\section{A soupçon of language} \label{sec:lang}
We begin by discussing some common language for last passage percolation models. We will abbreviate the phrase ``last passage percolation" as lpp. We consider lpp models in discrete time.

The purpose of this language is to understand similarities between formulas that appear in lpp models and tasep-like particle systems. We are thinking in particular about Sch\"{u}tz's formula for tasep \cite{Sch}, Johansson's formula for Geometric lpp \cite{Joh} and Warren's formula for Brownian tasep \cite{Warren}. These formulas are derived by different means. The motif, then, is to understand similarities between formulas by understanding similarities in the models themselves. For this, we need to understand what it means for models to be the same (isomorphism) and to be close (topology). Finally, we need to understand transformations between models and what is does to observables of interest.

\subsection{Domains for space, time and noise}
A lpp model comes with domains for space, time and noise. This is a triple
$$\cD = (D_s, D_t, D_{\eta})$$
defined as follows.

The temporal domain $D_t \subset \Z$ is a set of the form $D_t = [t_{min}, \infty)$ with the understanding that if $t_{min}= -\infty$ then $D_t = \Z$.
The domain $D_t$ inherits the linear order and discrete topology of $\Z$. ``Time" will be discrete throughout the discussion.

The spatial domain $D_s$ is a non-empty, closed subset of $\R$ together with its inherited linear order and topology.
Let $x_{min} = \inf \{y: y \in D_s\}$ with the understanding that $x_{min} = - \infty$ if $D_s$ is unbounded from below.
If $x_{min}$ is finite then it belongs to $D_s$ because $D_s$ is closed. Let $x^{-}_{min} = -\infty$. For every $x \in D_s$ with $x > x_{min}$, let
$$x^{-} = \sup \{y: y < x, y \in D_s\}.$$
The element $x^{-}$ is the predecessor of $x$. Note that $x^{-} \in D_s$ because $D_s$ is closed and the supremum is taken over a non-empty set that is bounded from above.
We assume that the mapping $x \in D_s\setminus \{x_{min}\} \mapsto x^{-} \in D_s$ is continuous.

The domain of noise $D_{\eta}$ is a non-empty, closed subset of $\R$ which is also a monoid, meaning that $0 \in D_{\eta}$ and if $x,y \in D_{\eta}$ then $x+y \in D_{\eta}$. The set $D_{\eta}$ inherits the linear order and topology from $\R$.

\subsection{Environment and Noise}
Given a domain $\cD$, a noise for lpp is a mapping
$$\eta : D_s \times D_t \to D_{\eta}$$
which is continuous. Here, $D_s \times D_t$ has the product topology. Thus, for every $t \in D_t$, $\eta(\cdot, t) : D_s \to D_{\eta}$ is a continuous function.
We set $\eta(x^-_{min},t) = 0$. We also assume that $\eta(x,t) - \eta(x^{-},t) \in D_{\eta}$ for every $x \in D_s$ and $t \in D_t$.

Noise is usually meant to be integrated over, and what we call noise is the integral of what it ought to be. The spatial increment $\eta(y,t) - \eta(x,t)$ represents the integral of the noise over the interval $(x,y]$.

The noise $\eta$ satisfies the no-negative-jumps condition if
$$ \eta(x,t) - \eta(x^-,t) \geq 0 \quad \text{for all}\;\; (x,t) \in D_s \times D_t.$$

A spatial interval $I$ is a subset $I = (x,y] \cap D_s = \{z \in D_s: x < z \leq y\}$. A temporal interval $J$ is a subset $J = [s,t] \cap D_t$.
An environment is a pair
$$ E = (I \times J, \eta)$$
where $I$ and $J$ are spatial and temporal intervals, respectively, and $\eta$ is the restriction of the noise over $I\times J$.
The set $I \times J$ is the support of the environment $E$.
We will abuse notation and write $(x,y]$ in place of $(x,y] \cap D_s$ and likewise for temporal intervals.

Two environments $E_1$ and $E_2$ are independent if their supports are disjoint: $(I_1 \times J_1) \cap (I_2 \times J_2) = \emptyset$.

Let $\cD$ be a domain. A noise $\eta : D_s\times D_t \to D_{\eta}$ is \textbf{stochastic} if there is a probability space $(\Omega, \mathcal{F}, \Pr)$ such that
$$\eta : (\Omega, \mathcal{F}, \Pr) \to \mathrm{C}(D_s\times D_t, D_{\eta})$$
is a random variable. Here $\mathrm{C}(D_s\times D_t, D_{\eta})$ is the space of continuous maps from $D_s\times D_t \to D_{\eta}$ with the topology of uniform convergence over compacts (which is metrizable).

A stochastic noise $\eta$ has the temporal independence property if the functions $\eta(\cdot, t)$ for $t \in D_t$ are mutually independent.
It has the spatial independence property if for every $t \in D_t$, and for every $x \leq y < u \leq v$ in $D_s$, the increment
$\eta(y,t)-\eta(x^-,t)$ is independent of the increment $\eta(v,t) - \eta(u^-,t)$.

\subsection{Directed paths}
Let $\cD$ be a domain. For $p = (x,s) \in D_s \times D_t$ and $q = (y,t) \in D_s \times D_t$, we write $p \leq q$ if $x \leq y$ and $s \leq t$.
Define
$$ D_{\uparrow}^2 = \{ (p;q) = (x,s;y,t) \in (D_s \times D_t)^2: p \leq q\}.$$

Suppose $p = (x,s)$ and $q = (y,t)$ with $p \leq q$. A path from $p$ to $q$ is a function
$$\pi: [s,t] \to D_s$$
such that $x \leq \pi(s) \leq \pi(s+1) \leq \cdots \leq \pi(t) = y$. By convention we set $\pi(s-1) = x$.
We say $\pi$ is monotone.

The graph of a path $\pi$ is the set
$$\Gamma(\pi) = \bigcup_{i=s}^t ([\pi(i-1)^{-}, \pi(i)]\cap D_s) \times \{i\} \subset D_s \times D_t$$
The graph is a compact set.

The collection of paths carries the Hausdorff topology on their graphs. This just means that $\pi_n \to \pi$ if there are times $s \leq t$ such that $\pi_n: [s,t] \to D_s$ is a monotone map for all large $n$, and $\pi_n(i) \to \pi(i)$ for every $i \in [s,t]$. For $p \leq q$, the collection of paths $\Pi_{(p;q)}$ from $p$ to $q$ is compact.

Let $\eta$ be a noise over the domain $\cD$. Let $p=(x,s) \leq q = (y,t)$ belong to $D_s\times D_t$.
The length of a path $\pi$ from $p$ to $q$ with respect to the noise $\eta$ is
$$ || \pi|| = \sum_{i=s}^t \eta(\pi(i),i) - \eta(\pi(i-1)^{-},i)$$
Note that $||\pi|| \in D_{\eta}$. The mapping  $ \pi \to ||\pi||$ is continuous.

\subsection{Last passage percolation}
An lpp model consists of the pair
$$ (\cD, \eta)$$
of domains for space, time, and noise together with the noise $\eta$ over $\cD$.
The last passage values are given by a function
$$L : D_{\uparrow}^2 \to D_{\eta}$$
defined according to 
$$L(p;q) = \max_{\pi \in \Pi_{(p;q)}} ||\pi||.$$
The maximum is achieved (though perhaps not uniquely) due to continuity of length and compactness of $\Pi_{(p;q)}$.
Note further that if $p = (x,s)$ and $q = (y,t)$ then $L(p;q)$ depends only on the environment $E = ((x^-,y] \times[s,t], \eta)$.
The function $L$ is continuous (here $D_{\uparrow}^2$ inherits its topology from $(D_s\times D_t)^2$).

One should think of $L$ as a directed metric, with $L(p;q)$ being the distance from $p$ to $q$. In this sense $L$ forms a length space: $L(p;q)$ is the length of the optimal path from $p$ to $q$. Although it may happen that $L(p;p) \neq 0$ if $\eta(x,t)-\eta(x^-,t) \neq 0$.

\subsection{Independent increments and composition law}
The function $L$ has two key properties. The first is independent increments along time. Suppose $S < T \in D_t$.
Consider the two processes
$$ \{ L(x,s;y,t): (x,s;y,t) \in D^2_{\uparrow}, s,t \leq S\} \quad \text{and} \quad \{ L(x,s;y,t): (x,s;y,t) \in D^2_{\uparrow}, s,t \geq T\}.$$
The former process depends on the environment $E_S = \big ( D_s \times [t_{min}, S], \eta \big)$ while the latter depends on $E_T = \big ( D_s \times [T, \infty), \eta \big)$. These two environments are independent. So we say that $L$ has independent increments along time.

The second property is the composition law. Suppose $s < r < t $ are times. Then for $x \leq y$,
\begin{equation} \label{eqn:complaw} L(x,s; y,t) = \max_{z \in [x,y] \cap D_s} \, \{L(x,s;z,r) + L(z,r;y,t) \}.\end{equation}
Indeed, simply condition on the location of optimal paths at time $r$.

\subsection{Isomorphisms of lpp}
\subsubsection{Maps between lpp models}
Consider two lpp models $(\cD^1,\eta_1)$ and $(\cD^2,\eta_2)$. A mapping $\Phi : (\cD^1,\eta_1) \to (\cD^2,\eta_2)$ is an object of the following kind.
It consists of a triple $\Phi = (\phi_1, \phi_2,\phi_3)$ with the following properties.
\begin{enumerate}
    \item The map $\phi_3: D^1_{\eta} \to D^2_{\eta}$ is continuous.
    \item The map $\phi_2: D^1_t \to D^2_t$ satisfies $\phi_2(s+1) = \phi_2(s)+1$. In other words, it is a translation.
    \item The map $\phi_1 : D^1_s \times \text{Im}(\phi_2) \to D^2_s$ is continuous. For every $t \in \text{Im}(\phi_2)$,
    the map $x \mapsto \phi_1(x,t)$ is surjective. For every $x, y \in D_s^1$ with $x \leq y$ and $t \in \text{Im}(\phi_2)$,
    $\phi_1(x,t) \leq \min \{ \phi_1(y,t), \phi_1(x,t+1)\}.$ Finally, $\phi(x^-,t) = \phi(x,t)^-$ for every $(x,t)$.
    \item Let $\Psi(x,t) = (\phi_1(x, \phi_2(t)), \phi_2(t)): D^1_s \times D^1_t \mapsto D^2_s \times D^2_t$.
    The following commutation relation must hold: $\phi_3(\eta_1(x,t)) = \eta_2(\Psi(x,t))$.
\end{enumerate}
For a mapping $\Phi$ denote by $\phi_3\cdot \eta_1$ the noise $\phi_3\cdot \eta_1(x,t) = \phi_3(\eta_1(x,t))$.
Denote by $\Phi \cdot L$ the lpp function for the model $(\cD^1, \phi_3\cdot \eta_1)$.

A mapping between two lpp models tries to capture a notion of "homomorphism". Note that the noises $\eta_1$ and $\eta_2$
are deterministic. Condition (3) above ensures that the map $\Psi$ takes environments to environments and directed paths to directed paths.
Condition (4) ensures that the noise $\eta_2$ over the image of $\Psi$ respects the noise $\eta_1$ up to a transformation.

\subsubsection{Composition of lpp maps}
Let $\Phi : (\cD^1,\eta_1) \to (\cD^2,\eta_2) $ and $\Phi': (\cD^2,\eta_2) \to (\cD^3,\eta_3)$ be maps between lpp models.
Their composition $\Phi \cdot \Phi'$ is a map from $(\cD^1,\eta_1) \to (\cD^3,\eta_3)$ as follows.
The map $\Phi \cdot \Phi' = (a,b,c)$ with
\begin{enumerate}
    \item $c = \phi'_3 \circ \phi_3$ and $b = \phi'_2 \circ \phi_2$.
    \item $a: D^1_s \times \text{Im}(b) \to D^3_s$ is as follows. Suppose $(x,t) \in D^1_s \times \text{Im}(b)$.
    Since $b$ is injective, there is a unique $t' \in D^1_t$ such that $b(t')=t$. So $t = \phi_2'(\phi_2(t'))$.
    Set $a(x,t) = \phi_1'(\phi_1(x,\phi_2(t')),t)$.
\end{enumerate}

It is not hard to verify that $\Phi'\cdot \Phi$ is also a valid mapping. The only non-trivial part is to check the commutation relation.
We have $c \cdot \eta_1(x,t) = \phi_3'(\phi_3(\eta_1(x,t))) = \phi_3'(\eta_2(\Psi(x,t))) = \eta_3(\Psi' \circ \Psi(x,t))$. Thus, we need to verify
$$ \Psi' \circ \Psi = \Psi_{\Phi' \cdot \Phi}.$$
Note that $a(x,b(t)) = \phi'_1(\phi_1(x,\phi_2(t)),b(t))$, and thus,
$$\Psi_{\Phi' \cdot \Phi}(x,t) = (\phi'_1(\phi_1(x,\phi_2(t)), \phi'_2 \circ \phi_2(t)), \phi'_2 \circ \phi_2(t)) = \Psi' \circ \Psi (x,t).$$

\subsubsection{Isomorphism}
Two lpp models $(\cD^1,\eta_1)$ ad $(\cD^2,\eta_2)$ are isomorphic if there are lpp maps $\Phi: (\cD^1,\eta_1) \to (\cD_2,\eta_2)$
and $\Phi^{-1}: (\cD^2,\eta_2) \to (\cD^1,\eta_2)$ such that $\Phi^{-1} \cdot \Phi = Id_{(\cD^1,\eta_1)}$ and $\Phi \cdot \Phi^{-1} = Id_{(\cD^2,\eta_2)}$.
Here $Id$ refers to the identity map of the corresponding lpp model. In this situation the maps $\Psi_{\Phi}$ and $\Psi_{\Phi^{-1}}$ are inverses of each other and we denote the latter by $\Psi^{-1}$.

Suppose $\pi: [s,t] \to D_s^1$ is a directed path. We find that $\phi_2([s,t]) = [\phi_2(s), \phi_2(s)+1,\ldots, \phi_2(s)+t-s]$.
For $0 \leq r \leq t-s$, let $u = \phi_2(s)+r = \phi_2(s+r)$. Define
$$ \Phi \cdot \pi(u) = \phi_1(\pi(s+r),u)$$
and note that $(\Phi \cdot \pi(u), u) = \Psi(\pi(s+r),s+r)$. Moreover,
$$ \Phi \cdot \pi(u+1) = \phi_1(\pi(s+r+1),u+1) \geq \phi_1(\pi(s+r+1),u) \geq \phi(\pi(s+r),u) = \Phi \cdot \pi(u).$$
If $\pi$ is a path from $p =(x,s)$ to $q = (y,t)$ then $\Phi \cdot \pi$ is a path from $\Psi(p)$ to $\Psi(q)$.
Due to the surjectivity of $x \mapsto \phi_1(x,t)$, any path from $\Psi(p)$ to $\Psi(q)$ is of the form $\Phi \cdot \pi$
for some path $\pi$ from $p$ to $q$. Also,
$$ || \Phi \cdot \pi||_{\eta_2} = \eta_2(\Psi(\pi(s),s) - \eta_2(\Psi(p^-,s)) + \cdots =
\phi_3 \cdot \eta_1(\pi(s),s) - \phi_3 \cdot \eta_1(p^-,s) + \cdots = || \pi||_{\phi_3 \cdot \eta_1}.$$
Consequently,
$$ \Phi \cdot L_{(\cD^1,\eta_1)}(p;q) = L_{(\cD^2,\eta_2)}(\Psi(p);\Psi(q)).$$
If $(\cD^1,\eta_1)$ is isomorphic to $(\cD_2,\eta_2)$ via $\Phi$ then the above gives
\begin{equation} \label{eqn:Isom} L_{(\cD^2,\eta_2)}(p;q) = \Phi \cdot L_{(\cD^1,\eta_1)}(\Psi^{-1}(p), \Psi^{-1}(q)).\end{equation}

If $(\cD^1, \eta_1)$ and $(\cD^2,\eta_2)$ are two lpp models with stochastic noises, then they are isomorphic if there exists a coupling $(\eta,\eta')$ of the two noises $\eta_1$ and $\eta_2$ such that $(\cD^1,\eta)$ is isomorphic to $(\cD^2,\eta')$ almost surely via a deterministic isomorphism $\Phi$.\footnote{This is not quite the ideal notion. It is better to consider the measures induced by the two stochastic noises over the space of lpp models, and then to define isomorphism with respect to these measures. But this gets tricky, so we are taking a way out by working with couplings.}

\subsection{Topology on lpp} \label{sec:topology}
Denote by $d_H$ the Hausdorff metric on compact subsets of $\R$. Let $C, C' \subset \R$ be closed and non-empty. Denote
$$ d_{set}(C, C') = \sum_{N=1}^{\infty} \frac{d_H(C \cap [-N,N], C'\cap[-N,N])}{2^N}.$$
For two compact subsets $K, K' \subset \R$, and $x \in K$, let $(x)_{K'} \in K'$ be the closest point in $K'$ to $x$.
Suppose $f: C \to \R$ and $f': C' \to \R$ are two continuous functions. For $N \geq 1$, let $K_N = C \cap [-N,N]$ and $K'_N = C' \cap [-N,N]$.
Set $$d_N(f,f') = \sup_{x \in K_N} |f(x) - f'((x)_{K'_N})| + \sup_{x \in K'_N} |f'(x) - f((x)_{K_N})|.$$
Define
$$d_{noise}(f, f') = \sum_{N=1}^{\infty} \frac{\min \, \{1, d_N(f,f')\}}{2^N}$$

Suppose $(\cD, \eta)$ and $(\cD',\eta')$ are two lpp models. Their distance is
$$ d((\cD, \eta); (\cD',\eta)) = | D_t \Delta D'_t| + d_{set}(D_s,D'_s) + d_{set}(D_{\eta}, D_{\eta'}) + \sum_{t \in \Z} \frac{d_{noise}(\eta(\cdot, t), \eta'(\cdot, t))}{2^{|t|}}$$
with the understanding that $\eta(t, \cdot) \equiv 0$ if $t \notin D_t$ and likewise of $\eta'$.
The function $d$ is a valid metric. (Although, it is not isomorphism invariant.)

A sequence of lpp models $(\cD_{\ell}, \eta_{\ell}) \to (\cD, \eta)$ as $\ell \to \infty$ if
$$ d((\cD_{\ell}, \eta_{\ell}); (\cD, \eta)) \to 0.$$

If $(\cD_1, \eta_1)$ and $(\cD_2,\eta_2)$ are two lpp models with stochastic noise, then the distance between them is the infimum of $\E{d((\cD_1,\eta); (\cD_2,\eta'))}$ where $(\eta,\eta')$ range over all couplings between the noises $\eta_1$ and $\eta_2$.
Let us denote this quantity $W_1((\cD_1, \eta_1); (\cD_2, \eta_2))$. This is the 1-Wasserstein distance.
A sequence of lpp models $(\cD_{\ell}, \eta_{\ell})$ with stochastic noises converge to an lpp model $(\cD, \eta)$ if $W_1((\cD_{\ell}, \eta_{\ell}); (\cD, \eta)) \to 0$.

\begin{prop} \label{prop:topo}
    Suppose a sequence of lpp models $(\cD_{\ell}, \eta_{\ell})$ converges to an lpp model $(\cD, \eta)$.
    If $(p_\ell,;q_\ell) \in D_{\ell, \uparrow}^2$ are such that $(p_\ell;q_\ell) \to (p;q) \in D_{\uparrow}^2$, then
    $L_{\ell}(p_\ell;q_\ell) \to L(p;q)$.
\end{prop}

\begin{proof}
    Let $p = (x,s)$, $q=(y,t)$, $p_n = (x_\ell,s_\ell)$ and $q_\ell = (y_\ell,t_\ell)$. Then $x_\ell \to x$, $y_\ell \to y$, and $s_\ell =s$ and $t_\ell=t$ for all large $\ell$. There is a compact interval $J = [-N,N] \subset \R$ such that for all large $\ell$, $L(p_\ell; q_\ell)$ depends on the environment
    $E_\ell = ((J\cap D_{s,\ell}) \times [s,t], \eta_{\ell})$ and $L(p;q)$ depends on $E = ((J\cap D_{s}) \times [s,t], \eta)$.
    
    Thus, $L_{\ell}(p_\ell; q_\ell) = || \pi_\ell||_{(\cD_{\ell}, \eta_{\ell})}$ where $\pi_\ell$ is some path from $p_\ell$ to $q_\ell$ whose graph lies inside $J \times [s,t]$. Likewise, $L(p;q) = ||\pi||_{(\cD,\eta)}$ where $\pi$ is a path from $p$ to $q$ whose graph  lies inside $J \times [s,t]$.
    The set of directed paths from $[s,t] \to J$ is compact. Any limit point of $\pi_{\ell}$ is a path $\pi$ from $p$ to $q$. Suppose $\pi'$ is such a limit point along the sequence $\pi_{\ell_k}$.
    Since length is continuous, and the distance $||\eta_{\ell} - \eta||$ between the noises are shrinking to zero, we have $||\pi_{\ell_k}||_{(\cD_{\ell}, \eta_{\ell})} \to ||\pi'||_{(\cD, \eta)}$. Also, $||\pi'||_{(\cD, \eta)} \leq L(p;q)$. Thus, $\limsup_{\ell} L_{\ell}(p_\ell;q_{\ell}) \leq L(p;q)$.

    Let $\pi$ be an optimal path from $p$ to $q$ for $(\cD, \eta)$ so that $||\pi||_{(\cD,\eta)} = L(p;q)$. We may find paths $\pi_{\ell}$ from $p_{\ell}$ to $q_{\ell}$ for $(\cD_{\ell}, \eta_{\ell})$ such that their graphs $\Gamma(\pi_{\ell}) \to \Gamma(\pi)$. Let $\delta_\ell = d((\cD_\ell, \eta_\ell); (\cD,\eta))$. There is a constant $C$ depending on $J \times I$ such that $| ||\pi_{\ell}||_{(\cD_\ell,\eta_\ell)} - ||\pi||_{(\cD,\eta)} | \leq C \delta_\ell$. Thus, $L_\ell(p_{\ell};q_{\ell} \geq ||\pi_{\ell}||_{(\cD_\ell,\eta_\ell)} \geq  ||\pi||_{(\cD,\eta)}|| - C\delta_{\ell} = L(p;q) - C\delta_{\ell}$. Therefore, $\liminf_{\ell} L_{\ell}(p_{\ell}; q_{\ell}) \geq L(p;q)$. 
\end{proof}

\begin{cor} \label{cor:weakconv}
    Let $(\cD_{\ell}, \eta_{\ell})$ be a sequence of lpp models with stochastic noises that converges to an lpp model $(\cD,\eta)$.
    For $1\leq k \leq n$, let $(p^k_\ell,;q^k_\ell) \in D_{\ell, \uparrow}^2$ be such that $(p^k_\ell;q^k_\ell) \to (p^k;q^k) \in D_{\uparrow}^2$.
    Then the $n$-tuples $(L_{\ell}(p^k_\ell;q^k_\ell); 1 \leq k \leq n) \to (L(p^k;q^k);1\leq k \leq n)$ in law.
\end{cor}

\subsection{Growth from a corner}
Let $(\cD, \eta)$ be an lpp model with last passage function $L$. Fix $(x_0,t_0) \in D_s \times D_t$.
Define the function
$$ G : ([x_0,\infty) \cap D_s) \times [t_0,\infty) \to D_{\eta}$$
according to 
\begin{equation} \label{eqn:G} G(x,t) = L(x_0,t_0; x,t).\end{equation}
Note that $G(x,t)$ depends on the environment supported over $(x_0^-, x] \times [t_0,t]$.
The function $G$ has both a temporal and spatial Markov property.

\subsubsection{Temporal Markov property}
Fix a time $T \geq t_0$. For every $t \geq T$, the process $G(\cdot,t)$ depends only on $G(\cdot, T)$ and the environment
$E = \big ((x_0^-,\infty] \times [T+1,t], \eta \big)$. The latter environment is independent of the environment
$E' = \big ((x_0^-,\infty] \times [t_0,T-1], \eta \big)$, upon which $G(\cdot, s)$ depends for every $s < T$.
This follows from the independent increments property of $L$.

In particular, if the noise $\eta$ is stochastic and has the temporal independence property, then $t \mapsto G(\cdot, t)$ is a Markov process along time.

\subsubsection{Spatial Markov property}
Fix a time $T > t_0$ or take $T = + \infty$. For $x \geq x_0$ in $D_s$, define
\begin{equation} \label{eqn:vecG} \vec{G}(x) = (G(x,t); t_0 \leq t < T) \in D_{\eta}^{[t_0,T)}.\end{equation}
The process $x \mapsto \vec{G}(x)$ is Markovian in the following sense.

Suppose $x_0 \leq x < y$ belong to $D_s$. Then for $t \geq t_0$,
$$G(y,t) = \max_{s \in [t_0,t]} \, \{G(x,s) + L(x,s; y,t) \}$$
This follows from the composition law by conditioning on the first time an optimal path from $(x_0,t_0)$ to $(y,t)$ hits spatial location $x$.
Observe that the quantities $L(x,s;y,t)$ above depend only on the environment supported over $(x^-,y]\times [t_0,t]$, which is independent
of the environment supported over $(x_0^-, x'] \times [t_0,t]$ upon which $G(x',s)$ depends for every $x' < x$ and $s \in [t_0,t]$.
We deduce that for $x < y$, $\vec{G}(y)$ depends only on $\vec{G}(x)$ and an environment independent of the environments that determine $\vec{G}(x')$ for every $x' < x$.
We call this the spatial Markov property of $G$.

If the noise $\eta$ is stochastic and has both the temporal and spatial independence properties, then $x \mapsto \vec{G}(x)$ is a Markov process.

\subsection{Space--time duality} \label{sec:spacetime}
Space-time duality transposes the role of space and time to build a particle system from an lpp model.

Let $(\cD, \eta)$ be an lpp model with corner growth function $G$ according to \eqref{eqn:G}. Assume the noise satisfies the no-negative-jumps condition:
$$\eta(x,t) - \eta(x^-,t) \geq 0 \quad \text{for every}\; (x,t) \in D_s \times D_t.$$
Due to continuity of $\eta$, the jumps of $\eta(x,t)$ come from gaps in $D_s$.

The space-time dual of $G$ is the process
\begin{equation} \label{eqn:Xdual} X_n(t) \end{equation}
defined for $n \in D_t$ with $n \geq t_0$ and $t \in D_s$ with $t \geq x_0$. It is simply
$$ X_n(t) = - G(t,n).$$
We think of $X$ as the evolution in time of a collection of particles indexed by integers. Colloquially,
$$X_n(t) = \text{position of particle number n at time t}.$$
The trajectory of particle number $n \geq t_0$ is the process $t \mapsto X_n(t)$.
The state of the particle system at time $t \in D_s$ with $t \geq x_0$ is the process $(X_n(t); n \geq t_0) \in D_{\eta}^{[t_0, \infty)}$.

Due to no-negative-jumps assumption,
$$ L(x_0,t_0, x, t) \geq L(x_0,t_0; x,t-1) + \eta(x,t)-\eta(x^-,t) \geq L(x_0,t_0; x,t-1).$$
Thus, $X_n(t) \leq X_{n-1}(t)$ for every $n > t_0$ in $D_t$ and $t \geq x_0$ in $D_s$.

By the composition law of lpp,
\begin{equation} \label{eqn:Grecursion}
    G(t,n) = \max_{r \in [x_0,t] \cap D_s} \, \{ G(r,n-1) + \eta(t,n) - \eta(r^-,n) \}
\end{equation}
which, then implies
\begin{equation} \label{eqn:Xrecursion} X_n(t) = \min_{r \in [x_0,t] \cap D_s}\,  \{X_{n-1}(r) +\eta(r^-,n) - \eta(t,n) \}.\end{equation}

The temporal Markov property of $G$ reduces to one-sided interaction of particles in $X$.
This means that the trajectory $X_n(\cdot)$ of particle number $n$ interacts only with the trajectories $X_k(\cdot)$ for $k < n$.
In particular, the trajectory of particle $n$ does not change if we remove particles $k > n$ from the system.
If the noise is stochastic and has the temporal independence property, then $X$ is Markovian in the ``spatial parameter" $n$.

The spatial Markov property of $G$ becomes a temporal Markov property for $X$.
This means that for every $ T \geq x_0$ in $D_s$, the state of the particle system at time $t > T$ depends only on the state at time $T$ and an environment that is independent of the environment that determines the state of the system at all times $s < T$.
If the noise is stochastic and has both the temporal and spatial independence properties, then $(X_n(t); n \geq t_0)$ is a Markov process in "time" $t \geq x_0$.

Under space-time duality, the ``events"
\begin{equation} \label{eqn:dualevents1}
    \{ X_{n_k}(t) \geq a_k; 1 \leq k \leq K \} = \{G(t,n_k) \leq -a_k; 1 \leq k \leq K \}.
\end{equation}
Therefore, if the noise is stochastic, the fixed time multi-space law of $X$ is equivalent to the fixed space multi-time law of $G$.

Likewise, the events
\begin{equation} \label{eqn:dualevents2}
    \{ X_{n}(t_k) \geq a_k; 1 \leq k \leq K \} = \{G(t_k,n) \leq -a_k; 1 \leq k \leq K \}.
\end{equation}
For stochastic noise, this translates the fixed particle multi-time law of $X$ to the fixed time multi-space law of $G$.

\subsection{Space--noise duality} \label{sec:spacenoise}
Space--noise duality builds a particle system out of an lpp model by inverting the roles of space and noise through a variational formula.
For this, one needs additional assumptions on the lpp model. This duality is what is frequently used to link tasep to lpp.

\subsubsection{Assumptions on the lpp model}
Let $(\cD, \eta)$ be an lpp model satisfying the following assumptions.
\begin{enumerate}
    \item $D_{\eta}$ is unbounded from above: $\sup \{y: y\in D_{\eta}\} = + \infty$.
    \item $D_s$ is bounded from below and unbounded from above: $x_{min} = \inf \{x: x \in D_s\}$ is finite and $\sup \{x: x\in D_{s}\} = + \infty$.
    By translating space, assume $x_{min} = 0$.
    \item The integer $0 \in D_t$ and for every $t \geq 0$, $\eta(x,t)$ satisfies: $\eta(0,t)=0$ and it is increasing in $x$ and tends to $+\infty$ as $x \to + \infty$.
\end{enumerate}

\subsubsection{Properties of the corner growth function}
Let $G$ denote the corner growth function \eqref{eqn:G} from $(x_0,t_0) = (0,0)$. It has the following properties due to the assumption made above.
\begin{enumerate}
    \item $G(0,n) = \eta(n,0) + \cdots + \eta(0,0) = 0$ for every $n \geq 0$.
    \item $G(y,n) \geq G(x,n) + \eta(y,n) - \eta(x^-,n) > G(n,x)$ for every $y > x \geq 0$ in $D_s$. Thus, $G(\cdot, n)$ is increasing.
    \item $G(x,n) \geq \eta(x,n) + \eta(0,n-1) + \cdots \eta(0,0) \geq \eta(x,n)$. Thus, $G(x,n) \geq 0$ and it is unbounded as $x \to +\infty$.
\end{enumerate}

\subsubsection{Particle system} \label{sec:particlesystem}
The space-noise dual is a particle system
$$ X^* = X^*_n(t)$$
defined for $n \geq 0$ in $D_t$ and $t \geq 0$ in $D_{\eta}$.
We think of $X^*_n(t)$ as the location of particle number $n$ at time $t$.

Firstly, let $X^*_n(0) \in \R$ for $n \geq 0$ be such that $$X^*_0(0) > X^*_1(0) > X^*_2(0) > \cdots.$$
We think of $X^*_n(0)$ as the initial position of particle number $n$. Then, let
\begin{equation} \label{eqn:Xstar}
    X^*_n(t) = X^*_n(0) + \sup \, \{ x \in D_s: G(x,n) \leq t\}.
\end{equation}
The process $X^*$ takes values in the set
$$ D^*_s = \{ x+y: x \in D_{s}, y = X^*_n(0)\;\text{for some}\; n\}.$$
We record some properties of $X^*$.

\begin{enumerate}
    \item At $t=0$, one has $G(x,n) \leq 0$ if and only if $x = 0$. Thus, $X^*_n(0) - X^*_n(0) = G(0,n) = 0$, as required.
    \item The set over which the supremum is taken in \eqref{eqn:Xstar} is non-empty for $t \geq 0$.
    It is bounded from above because $G(x,n)$ is unbounded as $x \to +\infty$. Since $x \mapsto G(x,n)$ is continuous, it is also closed.
    Thus, the supremum is attained within the set. Since $G(x,n)$ is also increasing, there exists an $x_n(t) \in D_s$ such that
    $$\{ x \in D_s: G(x,n) \leq t\} = [0,x_n(t)] \cap D_s.$$
    \item Since $G(x,n)$ is increasing in $x$, $x_n(s) \leq x_n(t)$ for $s \leq t$ in $D_{\eta}$. Thus, $t \mapsto X^*_n(t)$ is non-decreasing.
    \item By \eqref{eqn:Grecursion}, $G(x, n+1) \geq G(x, n) + \eta(x,n+1) - \eta(x^-,n) \geq G(x,n)$. Thus, if $t \geq G(x,n+1)$ then $t \geq G(x,n)$
    as well, which implies $X^*_{n+1}(t) - X^*_{n+1}(0) \leq X^*_n(t)-X^*_n(0)$. Since $X^*_{n+1}(0) < X^*_n(0)$, one has
    $$X^*_{n+1}(t) < X^*_n(t).$$
    \item The process $x_n(t)$ is determined by the environment supported over $D_s \times [0,n]$, and hence independent of the environment over $D_s \times [n+1,\infty)$.
    The temporal Markov property of $G$ implies $X^*_n(t)$ has one-sided interaction in the same sense as for the aforementioned space-time dual.
    \item The process $(X^*_n(t); n \geq 0) \in (D_s^*)^{[0,\infty)}$ is Markovian in time in the same sense as for the space-time dual. This is harder to prove and is not shown here.
    \item The ''events" $\{X^*_n(t) - X^*_n(0) \geq a\} = \{G(a,n) \leq t \}$. Therefore, for stochastic noise, the temporal law of $X^*_n(\cdot)$ is determined in terms of the spatial law of $G(\cdot, n)$.
    \item When $n=0$, if we assume $\eta(0,0) = 0$, then $X^*_0(t)-X^*_0(0) = \sup\, \{x \in D_s: \eta(x,0) \leq t\}$. This is the generalized inverse of $\eta(x,0)$: $\eta(X^*_0(t)-X_0(0),0) = t$.
\end{enumerate}

\section{Geometric lpp and discrete time tasep} \label{sec:geom}
Geometric last passage percolation is defined over the domain $D_t = D_s = D_{\eta} = \{0,1,2,3, \cdots\}$. The noise is semi-stochastic as it is given in terms of random weights $\omega_{i,j}$ for $i,j \geq 1$ together with a deterministic boundary condition. Consider independent random variables $\omega_{i,j}$ such that
$$ \Pr(\omega_{i,j} = k) = (1-q) q^{k-1} \quad k = 1,2,3,\cdots; \quad 0 < q < 1.$$
Consider a sequence of integers
$$ 0 = x_0 < x_1 < x_2 < x_3 < \cdots$$
Define the noise
$$ \eta(k, 0) = x_k $$
and, for $t > 0$,
$$ \eta(k, t) = \sum_{i=1}^k \omega_{i,t}, \quad \eta(0,t) = 0.$$
We are interested in the corner growth function
$$G(m,n) = L_{G}(0,0; m,n) \quad m,n \geq 0.$$
The function satisfies the recursion
\begin{equation} \label{eqn:lpprecur} G(m,n) = \max \{G(m-1,n), G(m,n-1)\} + \omega_{m,n}\end{equation}
with boundary conditions $G(0,n) = 0$ and $G(m,0) = x_m$.
Since $x_k \geq 0$, $\omega_{i,j} \geq 1$ and $\eta(0,\cdot) = 0$, an optimal path from $(0,0)$ to $(m,n)$ will always start by taking a step rightward instead of upward.
Thus,
\begin{equation} \label{eqn:geompolymer} G(m,n) = \max_{1 \leq k \leq m} \, x_k + L_G(k,1; m,n).\end{equation}

The utility of Geometric lpp is two-fold. Firstly, the model has strong symmetries due to space, noise and time occupying the same domain.
Secondly, the model allows for limit transitions to several other models of lpp, such as Exponential lpp, Brownian lpp, Poissonian lpp, and the directed landscape.

\subsection{Johansson's formula}
Consider a transformed lpp on the same domain but with the noise given by random weights
$$\omega'_{i,j} = \omega_{i,j} - 1 \quad \text{for}\; i,j \geq 1$$
and $\eta'(0,t) = \eta'(k,0) \equiv 0$.
Let $G'$ denote the growth function from $(0,0)$. It satisfies
$$G'(m,n) = \max \{G'(m-1,n), G'(m,n-1)\} + \omega'_{m,n}$$
with boundary conditions $G'(0,n) = 0$ and $G'(m,0) = 0$.
Johansson \cite{Joh} found a formula for the transition probabilities of the Markovian evolution of $G'$ along space. It goes as follows.

Let $N \geq 1$ be large. Define the random vector
$$ \vec{G}'_s(m) = (G'(m,1), G'(m,2), \ldots, G'(m,N)).$$
It takes values in the set
$$ \mathbb{W}^N = \{(x_1, \ldots, x_N) \in \mathbb{Z}^N: 0 \leq x_1 \leq x_2 \leq \cdots \leq x_N\}.$$
For $x,y \in \mathbb{W}^N$ and $m > \ell$, Johansson proved:
\begin{equation} \label{eqn:transition}
P(\vec{G}'_s(m) = y \mid \vec{G}'_s(\ell)=x) = \det [\nabla^{j-i}w_{m-\ell}(y_j-x_i)].
\end{equation}
Here $\nabla$ is the discrete differential: $\nabla f(x) = f(x+1) - f(x)$ with inverse $\nabla^{-1}f(x) = \sum_{y < x} f(y)$.
The function $w_n(x)$ is the $n$-step transition probability of a random walk with $Geom(1-q)-1$ step distribution:
$$ w_n(x) = P( X_1 + \cdots + X_n = x), \quad X_i \;\text{are i.i.d.}\; Geom(1-q)-1.$$
There is an explicit formula (Negative Binomial formula):
$$ w_n(x) = \binom{x+n-1}{n-1}q^x(1-q)^n \mathbf{1}_{x \geq 0}.$$

\subsection{Sch\"{u}tz's formula}
Consider the continuous time tasep particle system on $\Z$ with $N$ particles. Denote by $X_1(t) > X_2(t) > \cdots > X_N(t)$ the location of the particles at time $t$.
Sch\"{u}tz \cite{Sch} found a formula for the transition probability of $(X_1(t), \ldots, X_N(t))$. It reads as follows.
\begin{equation} \label{eqn:schutz}
    \Pr(X_1(t) = y_1, \ldots, X_N(t) = y_N \mid X_1(0) = x_1, \ldots, X_N(0) = x_N) = \dt{F_{i-j}(y_{N+1-i}-x_{N+1-j})}
\end{equation}
Here $y_1 > \cdots > y_N$ and $x_1 > \cdots > x_N$ are integers. The function $F_k(x)$ is
\begin{equation} \label{eqn:FSch} F_k(x) = \frac{(-1)^k}{2\pi \mathbf{i}} \oint_{|w|=r>1} \frac{dw}{w}\, \frac{w^{k-x}}{(1-w)^k}\, e^{(w-1)t} \end{equation}

This formula has similarities with Johansson's formula \eqref{eqn:transition}.
It was used by \cite{BFPS, Sas} to write the probability
\begin{equation} \label{eqn:tasepcorrelation} \Pr(X_{n_1}(t) > a_1, \ldots, X_{n_k}(t) > a_k)\end{equation}
as a Fredholm determinant by using the Eynard-Mehta method. The authors came to realize the process $(X_1(t), \ldots, X_N(t))$
as the edge of a random Gelfand-Tseltin pattern, the law of which is a 2-d determinantal process.
Then, there is an correlation kernel $K$ (coming from the Eynard-Mehta method) such that \eqref{eqn:tasepcorrelation}
is the Fredholm determinant of $K$ because \eqref{eqn:tasepcorrelation} is a gap probability of the aforementioned determinantal process.
In order to get a explicit formula for $K$ one has to solve a biorthogonalization problem.
The biorthogonalization problem depends on the initial condition of $X$, and such a problem was solved for tasep in \cite{MQR},
enabling one to express \eqref{eqn:tasepcorrelation} as an explicit Fredholm determinant for all initial conditions of tasep.

\subsection{Markovian transition of $G$}
We wish to find the transition probability of the temporal evolution of $G$ and use it to derive the probability
\begin{equation} \label{eqn:Gtemporal}
    \Pr(G(m_1,n) < a_1, \ldots, G(m_k,n) < a_k)
\end{equation}
as a Fredholm determinant.

Looking at the growth function $G'$, observe that the transposition of space-time: $(m,n) \mapsto (n,m)$,
is an isomorphism for the corresponding lpp model.
Here, the random weights $\omega'_{i,j}$ are coupled to themselves by placing weight $ \omega'_{j,i}$ at location $(i,j)$, which preserves their joint law.
Thus, the formula \eqref{eqn:transition} also serves as the transition probability of the temporal law of $G'$.
With $N \geq 1$ large, define
$$ \vec{G}'_t(n) = (G'(1,n), \ldots, G'(N, n)).$$
Then for $x,y \in \mathbb{W}^N$ and $n > \ell$:
\begin{equation} \label{eqn:transition2}
\Pr(\vec{G}'_t(n) = y \mid \vec{G}'_t(\ell)=x) = \det [\nabla^{j-i}w_{n-\ell}(y_j-x_i)].
\end{equation}

The relation between $G'$ and $G$ is $G(m,n) = G'(m,n) + m +n\ind{m\neq 0}$.
We have imposed the boundary condition $G(k,0) = x_k$, which become the boundary conditions $G'(k,0)=x_k-k$.
Let $$\vec{G}_t(n) = (G(1,n), \ldots, G(N, n)).$$
Observe that $\vec{G}_t(0)= (x_1, x_2, \ldots, x_N)$.
Then for $y = (y_1, \ldots, y_N)\in \Z^N$ with $0 \leq y_1 < y_2 < \cdots < y_N$, \eqref{eqn:transition2} implies
\begin{equation} \label{eqn:transition3}
    \Pr(\vec{G}_t(n) = y) = \det [\nabla^{j-i}w_n(y_j-x_i +i-j -n)].
\end{equation}

The next step is to represent the right side of \eqref{eqn:transition3} as a Sch\"{u}tz-like formula \eqref{eqn:schutz}.
This would allow one to use the methods of \cite{BFPS, Sas, MQR} to express \eqref{eqn:Gtemporal} as a Fredholm determinant.
The idea is as follows. Let $\vec{G}_s(m) = (G(m,1), \ldots, G(m,N))$ be the spatial Markov process of $G$.
Apply the space-time transposal isomorphism so that the boundary condition goes from being at the bottom to the left.
Thus, $\vec{G}_s(m)$ evolves from $\vec{G}_s(0) =x$. Consider the space-time dual $X$ of $G$ from \eqref{eqn:Xdual}.
One has $X(t) = (X_1(t), \ldots, X_N(t)) = - (G(t,1), \ldots, G(t,N)) = - \vec{G}_s(t)$. Now,
\begin{equation} \label{eqn:XGrelation}
\Pr(\vec{G}_t(n) = y \mid \vec{G}_t(0) = x) = \Pr(\vec{G}_s(n) = y \mid \vec{G}_s(0)=x) = \Pr(X(n) = -y\mid X(0) = -x)
\end{equation}
We can thus write the transition probability of $X$ using \eqref{eqn:transition3}.
If we express it in terms of the number-reversed process $X_{N+1-m}(n)$, $1\leq m\leq N$,
then the resulting formula has the desired form.

\begin{prop} \label{prop:JoSc}
    Define, for $k \in \Z$ and $x \in \Z$, the function
    \begin{equation}
    F_k(x) = \frac{(-1)^k}{2\pi \mathbf{i}} \oint_{|w|=r>1} \frac{dw}{w}\, \frac{w^{k-x}}{(1-w)^k}\, \left (\frac{(1-q)}{w-q} \right)^n
\end{equation}
Referring to \eqref{eqn:transition3}, we have
$$\Pr(\vec{G}_t(n)=y) = \det [F_{i-j}(\tilde{y}_{N+1-i}-\tilde{x}_{N+1-j})]$$
where $\tilde{y}_k = -y_k$ and $\tilde{x}_k = -x_k$.
\end{prop}

\begin{proof}
    We execute the idea above with a direct computation. Let $H_k(x) = \nabla^k w_n(x)$. Using the generating function of negative Binomial coefficients, one can write $H_k(x)$ as a contour integral (see \cite{Joh}).
    \begin{equation} \label{eqn:Hk}
   H_k(x) = \frac{(-1)^{k-1}}{2 \pi \mathbf{i}} \oint_{|z| = r > 1} dz\, \frac{z^k (1-z)^{n+x-1}}{(1 - \frac{z}{1-q})^n}.
\end{equation}
One has $F_k(x) = H_{-k}(k-n-x)$. This follows from \eqref{eqn:Hk} by changing variables $z \mapsto 1-w$ in the integral.

Define $\tilde{F}_k(x) = H_{-k}(-x)$. Then $\tilde{F}_k(x-k+n) = F_k(x)$. Now we compute, by \eqref{eqn:transition3}:
\begin{align*}
  \Pr(\vec{G}_t(n) = y) &= \det [H_{j-i}(y_j-x_i +i-j-n)] \\
  &= \det [\tilde{F}_{i-j}(x_i-y_j+j-i+n)] \\
  &= \det [\tilde{F}_{j-i}(x_{N+1-i}-y_{N+1-j}+i-j+n)]\\
  &= \det [\tilde{F}_{i-j}(x_{N+1-j}-y_{N+1-i} +j-i+n)]\\
  &= \det [\tilde{F}_{i-j}(\tilde{y}_{N+1-i} - \tilde{x}_{N+1-j} +j-i+n)]\\
  &= \det [F_{i-j}(\tilde{y}_{N+1-i} - \tilde{x}_{N+1-j})]
\end{align*}
\end{proof}

\subsection{Fredholm determinants and kernels} \label{sec:fredholm}
Let $K$ be an integral kernel acting on the space $L^2(\Omega,\mu)$. The Fredholm determinant of $K$ is
\begin{equation} \label{eqn:Fredholm}
\dt{I+K} = 1 + \sum_{k=1}^{\infty} \frac{1}{k!} \int_{\Omega^K} d\mu(z_1) \cdots d\mu(z_k)\, \dt{K(z_i,z_j)}_{1\leq i,j\leq k.}\end{equation}
If there are functions $f$ and $g$ on $X$ such that
$$ |K(u,v)| \leq f(u) g(v)$$
with $f$ bounded and $g$ being integrable (or vice versa), then the series converges absolutely. Furthermore, suppose a sequence of integral kernels
$K_n$ satisfy $K_n \to K$ pointwise on $\Omega$ and $|K_n(u,v)| \leq f(u) g(v)$ for every $n$. Then
$$ \dt{I+K_n}_{L^2(\Omega,\mu)} \to \dt{I+K}_{L^2(\Omega,\mu)}$$
See \cite{Josurvey, JoDPG} for proofs of these facts, which are deduced from Hadamard's inequality and dominated convergence theorem.

\subsubsection{Matetski's and Remenik's kernel}
Consider the space $\Omega = \{1,2, \ldots, k\} \times \Z$ with the counting measure $\mu$. Define an integral kernel $K^*_N$ on $\Omega$ as follows.
Let $0 = x_0 < x_1 < x_2 < \cdots$ be increasing integers. Let $\tilde{x}_i = -x_i$. Let $1 \leq m_1 , m_2 , \cdots , m_k \leq N$ and $n \geq 1$ be integers. Assume the $m_i$ are distinct.
Choose $\theta \in (0,1)$ and set $\alpha = (1-\theta)/\theta$.
Define
$$ \phi(w) = \frac{1-q}{w-q}; \quad 0 < q < 1 \;\text{and}\; w \in \C\setminus{\{q\}}.$$
Let $r \in (q,1)$ and $\delta < 1$ be radii parameters.

Define the kernel
$$Q^*(z_1,z_2) = (1-\theta) \theta^{z_1-z_2-1} \mathbf{1}_{z_1 > z_2}$$
This is the transition matrix of a random walk with $-Geom(1-\theta)$ steps [steps strictly to the left]. Let $B^*_m$ denote the corresponding random walk.
For integers $m \in \mathbb{Z}$ and $n \geq 1$, define the following kernels.
\begin{equation}
(Q^*)^m(z_1,z_2) = \frac{1}{2\pi \mathbf{i}} \oint_{|w|=r} dw\, \frac{\theta^{z_1-z_2}}{w^{z_1-z_2-m+1}} \left(\frac{\alpha}{1-w}\right)^m.
\end{equation}

\begin{equation}S_{n,-m}^*(z_1,z_2) = \frac{1}{2 \pi \mathbf{i}} \oint_{|w|=r}dw\, \frac{\theta^{z_1-z_2}}{w^{z_1-z_2+m+1}} \left( \frac{1-w}{\alpha}\right)^m \phi(w)^n.\end{equation}

\begin{equation}\bar{S}^*_{n,m}(z_1,z_2) = \frac{1}{2 \pi \mathbf{i}} \oint_{|w|=\delta}dw\, \frac{\theta^{z_1-z_2} (1-w)^{z_2-z_1+m-1}}{(w/\alpha)^m} \phi(1-w)^{-n}.\end{equation}
Let
$$\tau_N = \min \, \{m=0,1,\ldots,N-1: B^*_m > \tilde{x}_{m+1}\}$$
with $\tau_N = +\infty$ if $B^*$ does not hit $\tilde{x}$ by time $N-1$.
Define the kernel
\begin{equation}S^{epi(\tilde{x})}_{n,m}(z_1,z_2) = \E{ \bar{S}^*_{n,m-\tau_N}(B^*_{\tau_N},z_2) \mathbf{1}_{\tau_N < m} \mid B^*_0 = z_1} \end{equation}
Finally, define the kernel
\begin{equation} \label{eqn:Kgeometric}
K^*_N(i,\cdot; j, \cdot) = - (Q^*)^{m_j-m_i}\mathbf{1}_{m_i < m_j} + S_{n,-m_i}^{*} \cdot S^{epi(\tilde{x})}_{n,m_j}
\end{equation}

These kernels were introduced by Matetski and Remenik \cite{MR} in order to study transition probabilities for tasep-like particle systems.
The have a certain structure.
For $n \in \Z$, define the kernels
$$R_{n}(z_1, z_2) = \frac{1}{2 \pi \mathbf{i}} \oint_{|w|=r}dw\, \frac{\theta^{z_1-z_2}}{w^{z_1-z_2+1}} \phi(w)^{n}.$$
The $R_n$s form an abelian group: $R_{n+m} = R_n \cdot R_m$ and $R_0 = Id$.
One has the relations:
\begin{align*}
    S^*_{n,-m} &= Q^{-m}R_n \\
    \bar{S^*}_{n,m} &= \bar{Q}^{(m)} R_{-n}, \quad \bar{Q}^{(m)} = \bar{S^*}_{0,m}.
\end{align*}

\subsubsection{Kernel for Geometric lpp with boundary conditions}
Define the following kernel $K$ on $\Omega$ using the notation from the previous section.
\begin{equation} \label{eqn:Kgeom}
    K(i,\cdot; j, \cdot) = - Q^{m_j-m_i}\mathbf{1}_{m_i < m_j} + S_{n,-m_i} \cdot S^{hypo(x)}_{n,m_j}
\end{equation}
The kernels $Q$, $S$ and $S^{hypo(x)}$ are as follows for $m,n \in \Z$ with $n \geq 1$.
\begin{equation} \label{eqn:Qgeom}
Q^m(z_1,z_2) = \frac{1}{2\pi \mathbf{i}} \oint_{|w|=r < 1} dw\, \frac{\theta^{z_2-z_1}}{w^{z_2-z_1-m+1}} \left(\frac{\alpha}{1-w}\right)^m.
\end{equation}
This is the $m$-step transition probability of a random walk $B_m$ with $Geom(1-\theta)$ steps [strictly to the right].
\begin{equation} \label{eqn:Sgeom}
S_{n,-m}(z_1,z_2) = \frac{1}{2 \pi \mathbf{i}} \oint_{|w|=r > q}dw\, \frac{\theta^{z_2-z_1}}{w^{z_2-z_1+m+1}} \left( \frac{1-w}{\alpha}\right)^m \phi(w)^n.
\end{equation}
\begin{equation} \label{eqn:Sbargeom}
\bar{S}_{n,m}(z_1,z_2) = \frac{1}{2 \pi \mathbf{i}} \oint_{|w|=\delta < 1}dw\, \frac{\theta^{z_2-z_1} (1-w)^{z_1-z_2+m-1}}{(w/\alpha)^m} \phi(1-w)^{-n}.
\end{equation}
Let
$$\tau = \min \, \{m=0,1,2, \ldots: B_m < x_{m+1}\}$$
with $\tau = +\infty$ if $B$ does not hit $x$.
The kernel
\begin{equation} \label{eqn:Shypogeom}
S^{hypo(x)}_{n,m}(z_1,z_2) = \E{ \bar{S}_{n,m-\tau}(B_{\tau},z_2) \mathbf{1}_{\tau < m} \mid B_0 = z_1}.
\end{equation}

\subsection{Fixed time law of Geometric lpp with boundary conditions}
\begin{thm} \label{thm:geomlpp}
    Consider the Geometric lpp model \eqref{eqn:geompolymer} with integer boundary condition $0 < x_1 <x_2 < x_3 < \cdots$.
    Let $n \geq 1$, $m_1, \ldots, m_k \geq 1$ and $a_1, \ldots, a_k$ be integers such that the $m_i$ are distinct. Then,
    \begin{equation} \label{eqn:geomlpp} \Pr(G(m_1,n) < a_1, \ldots, G(m_k,n) < a_k) = \dt{I - \chi_a K \chi_a}_{\ell^2(\{1,2,\ldots,k\}\times \Z)}.\end{equation}
    The kernel $K$ is from \eqref{eqn:Kgeom} and $\chi_a(i,z) = \ind{z \geq a_i}$.
\end{thm}

\begin{proof}
    This is an application of Theorem 1.2 of \cite{MR}. Fix $N > \max \{m_1, \ldots, m_k\}$ and consider $G(n) = (G(1,n), \ldots, G(N,n)$.
    The probability in \eqref{eqn:geomlpp} is a marginal of $G(n)$. It also does not depend on $N$ when $N > \max \{m_1, \ldots, m_k\}$.
    By \eqref{eqn:XGrelation} and Proposition \ref{prop:JoSc},
    $$ \Pr(G(n)=y\mid G(0)=x) = \Pr(X(n)=\tilde{y}\mid X(0)=\tilde{x}) = \dt{F_{i-j}(\tilde{y}_{N+1-i} - \tilde{x}_{N+1-j})}.$$
    By Theorem 1.2 of \cite{MR},
    $$\Pr(X_{m_1}(n) > -a_1, \ldots, X_{m_k}(n) > - a_k) = \dt{I - \bar{\chi} K_N^{*} \bar{\chi}}_{\ell^2(\{1, \ldots, k\} \times \Z)}.$$
    Here $\bar{\chi}(i,z) = \ind{z \leq -a_i}$ and $K^*_N$ is from \eqref{eqn:Kgeometric}.
    The probability on the left is the probability in \eqref{eqn:geomlpp}.

    Since the probability does not depend on $N$, consider the limit $N \to \infty$ of the Fredholm determinant.
    The only way $K^*_N$ depends on $N$ is through the hitting time $\tau_N = \min \{ m=0,1,\ldots,N-1: B^*_m > \tilde{x}_{m+1} \}$.
    Define $\tau^* = \inf \{ m \geq 0: B^*_m > \tilde{x}_{m+1}\}$. We can couple all the $\tau_N$ and $\tau^*$ together by using
    the same random walk $B^*$. Under this coupling, $\tau_N \downarrow \tau^*$ almost surely. Moreover, for any fixed $m > 0$,
    $\ind{\tau_N < m}$ stabilizes to $\ind{\tau^* < m}$ as $N \to \infty$. Thus, $K^*_N$ stabilizes to the kernel $K^*$
    for all large $N$, where $K^*$ has the same definition as $K^*_N$ but with $\tau^*$ in place of $\tau_N$.
    Consequently, $\dt{I - \bar{\chi} K_N^{*} \bar{\chi}}$ equals $\dt{I - \bar{\chi} K^{*} \bar{\chi}}$ for all sufficiently large $N$.

    Finally, change variables $z \mapsto -z$ in the integrals of the Fredholm determinant expression \eqref{eqn:Fredholm}
    for $\bar{\chi}K^*\bar{\chi}$. This changes to the Fredholm determinant expression for $\chi_a K \chi_a$.   
\end{proof}

\subsection{Discrete time tasep}
Discrete time parallel update tasep is an interacting particle system that evolves as follows.
At time $t=0$ there is some initial configuration of particles on the integer lattice $\Z$.
At each time $t > 0$, a particle attempts to jump one unit rightward with probability $1-q$, independently of others.
The jump is successful if there is no particle at that location at time $t-1$.

Let us focus on the case when there is a rightmost particle, which we call particle number $0$.
Label the particles from right to left by integers $n = 0,1,2,3, \cdots$. Let $Y_n(t)$ be the location
of particle number $n$ at time $t$. Let
$$ G_Y(m,n) = \text{time when particle number}\; n\; \text{makes jump number}\; m.$$
for $n \geq 0$ and $m \geq 1$. It is clear that $G_Y$ together with the initial condition $Y_n(0)$
determine the evolution of tasep. The function $G_Y$ satisfies the recursion
\begin{equation} \label{eqn:dtasep}
G_Y(m,n) = \max \{G_Y(m-1,n), G_Y(m+1-\mathrm{gap}_n, n-1) \} + \omega_{m,n}
\end{equation}
where $\omega_{m,n}$ are independent $Geom(1-q)$ random variables and $\rm{gap}_n = Y_{n-1}(0) - Y_n(0) \geq 1$.

The case $\rm{gap}_n \equiv 1$ is called the step initial condition. Assuming $Y_0(0)=0$, this means $Y_n(0) = -n$ for $n \geq 0$.
The recursion \eqref{eqn:dtasep} becomes \eqref{eqn:lpprecur} and we find that $G(m,n)$ from $\eqref{eqn:lpprecur}$ is the time
when particle $n$ makes its $m$-th jump. The process $Y$ is the space-noise dual $X^*$ of $G$ according to \eqref{eqn:Xstar}.
For $n=0$ we find that $Y_0(t)$ is a deterministic trajectory with particle $0$ jumping one unit rightward at the times $0 < x_1 < x_2 < x_3< \cdots$,
where $x$ is the boundary condition for $G$ in \eqref{eqn:lpprecur}. Theorem \ref{thm:geomlpp} then has the following corollary.

\begin{cor} \label{cor:dtasep}
    Consider discrete time parallel update tasep from the step initial condition $X^*_n(0) = -n$ for $n \geq 0$.
    Let $X^*_n(t)$ be the location of particle number $n$ at time $t$. Let $X^*_0$ jump one unit rightward at the integer times $x_k$ with $0 < x_1 < x_2 < x_3 < \cdots$. The other particles move according to the stochastic rule described above. Then for $n \geq 1$, integer times $t_1, t_2, \ldots, t_k \geq 1$, and distinct integers $a_1, \ldots, a_k \geq -n$,
    \begin{equation} \label{eqn:dtasepformula}
    \Pr(X^*_n(t_i) > a_i; 1 \leq i \leq k) = \dt{I - \chi_t J \chi_t}_{\ell^2(\{1,2,\ldots,k\}\times \Z)}
    \end{equation}
    where $\chi_t(i,z) = \ind{z > t_i}$ and the kernel $J$ is the kernel $K$ from Theorem \ref{thm:geomlpp} with choice of parameters $m_i = a_i+n+1$ and $n$ the same.
\end{cor}

\begin{proof}
    Recall property (7) of the space-noise dual $X^*$ from Section \ref{sec:particlesystem}. The event $\{X^*_n(t) > a\}$ is equivalent to the event $\{G(a+n+1,n) < t+1\}$. The corollary follows from Theorem \ref{thm:geomlpp}.
\end{proof}

\subsection{A word about space-time and space-noise duality}
Consider again the tasep-like process $X$ from \eqref{eqn:Xdual} that is the space-time dual of Geometric lpp.
By \eqref{eqn:XGrelation} and Proposition \ref{prop:JoSc}, its transition probabilities are
\begin{equation}
    \Pr(X(t) = y \mid X(0)=x) = \dt{F_{i-j}(y_{N+1-j} - x_{N+1-i})}.
\end{equation}
The recursion \eqref{eqn:Grecursion} implies that $X$ satisfies the recursion
$$ X_n(t) = \min \, \{X_{n-1}(t), X_n(t-1)\} - \omega_{n,t}$$
where the $\omega_{n,t}$ are independent $Geom(1-q)$ random variables.
This process is called the geometric PushTASEP with left jumps (see Section 2.3 of \cite{MR} and Case A of \cite{DW}).
The particles move leftward and push the particles in front of them to maintain their ordering.

Now consider the space-noise dual $X^*$ of $G$. Here the particles move rightward and are blocked by the particles in front of them. For the step initial condition, the recursion satisfied by $X^*_n(t)$ (discrete time parallel update tasep) is
$$ X_n(t) = \min \, \{X_n(t-1) + \xi_{t,n}, X_{n-1}(t-1)-1\} \quad \xi_{t,n} \;\text{are i.i.d.}\; Bernoulli(1-q).$$
Duality, in a concrete sense, links blocking to pushing.

\section{Exponential lpp and continuous time tasep} \label{sec:exp}
Exponential last passage percolation is defined over the domain $D_t = D_s = \{0,1,2,3, \cdots\}$ and $D_\eta = [0,\infty)$ being all non-negative real numbers. The noise is again semi-stochastic and given in terms of random weights $\omega_{i,j}$ for $i,j \geq 1$ together with a boundary condition. Consider independent Exponential random variables $\omega_{i,j}$ such that
$$ \Pr(\omega_{i,j} > x) = e^{-x}, \quad x \geq 0.$$
Consider a sequence of real numbers
$$ 0 = x_0 < x_1 < x_2 < x_3 < \cdots \to +\infty.$$
Define the noise
\begin{align} \label{eqn:Enoise}
    \eta(k, 0) & = x_k \quad \text{for}\; k \geq 0\\
    \eta(k, t) &= \sum_{i=1}^k \omega_{i,t}, \quad \eta(0,t) = 0 \quad \text{for}\; t > 0, k \geq 0.
\end{align}

We are again interested in the corner growth function
$$G(m,n) = L_{E}(0,0; m,n) \quad m,n \geq 0.$$
The function satisfies the recursion
\begin{equation} \label{eqn:Elpprecur} G(m,n) = \max \{G(m-1,n), G(m,n-1)\} + \omega_{m,n}\end{equation}
with boundary conditions $G(0,n) = 0$ and $G(m,0) = x_m$.
Since $x_k > 0$ and $\omega_{i,j} \geq 0$, an optimal path from $(0,0)$ to $(m,n)$ will always start by taking a step rightward instead of upward.
Thus, as before,
\begin{equation} \label{eqn:exppolymer} G(m,n) = \max_{1 \leq k \leq m} \, x_k + L_E(k,1; m,n).\end{equation}

For later convenience, we shall condition the model on the almost sure event that $\{\omega_{i,t} > 0 \;\text{for all}\; (i,t)\}$.

\subsection{Geometric to Exponential transition}
Exponential lpp is a limit of Geometric lpp as the parameter $q \to 1$. Suppose $\omega_q \sim Geom(1-q)$.
Then $(1-q)\omega_q \to \omega$ in law as $q\to 1$ with $\omega \sim Exp(1)$.
Let $\ell \geq 2$ be an integer and set $q_{\ell} = 1 - \frac{1}{\ell}$. Set $x_k^{\ell} = \lfloor \ell x_k \rfloor / \ell$ and note that
$\sup_k \, |x_k^{\ell} - x_k| \leq 1/\ell$. Let $\tilde{\omega}_{i,j}^{\ell}$ be i.i.d. $Geom(1-q_{\ell})$ random variables. Set $\omega_{i,j}^{\ell} = \tilde{\omega}_{i,j}^{\ell} / \ell$.

Consider the lpp model on domain $\cD^{\ell}$ where $D^{\ell}_s = D^{\ell}_t = \{0,1,2,\ldots\}$, and
$D^{\ell}_{\eta} = \{0, \ell^{-1}, 2\ell^{-1}, \ldots\}$. The noise $\eta^{\ell}$ is defined according to \eqref{eqn:Enoise} using the weights $\omega^{\ell}_{i,j}$ and $x_k^{\ell}$. There is a coupling of the random variables $\tilde{\omega}_{i,j}^{\ell}$ for all $i,j \geq 1$ and $\ell \geq 2$
together with i.i.d. random variables $\omega_{i,j} \sim Exp(1)$ such that $\omega^{\ell}_{i,j} \to \omega_{i,j}$ almost surely as $\ell \to \infty$,
simultaneously for all $i,j$. This follows from Skorokhod representation theorem.
As a result, the lpp model $(\cD^{\ell}, \eta^{\ell})$ converges to the Exponential lpp model $(\cD, \eta)$ above as $\ell \to \infty$ in the sense of Section \ref{sec:topology}.

Due to the convergence of models, one has the convergence of the probabilities:
\begin{equation} \label{eqn:geomtoE}
\Pr(G(m_1,n) \leq a_1, \ldots, G(m_k,n) \leq a_k) = \lim_{\ell} \; \Pr(G^{\ell}(m_1,n) \leq a_1, \ldots, G^{\ell}(m_k,n) \leq a_k)\end{equation}
where $G^{\ell}$ is the corner growth function of $(\cD^{\ell},\eta^{\ell})$. This follows from Corollary \ref{cor:weakconv}.

In fact, one can also show convergence of the Markovian transition probabilities using the formula \eqref{eqn:transition}.
Let $G(n) = (G(1,n), \ldots, G(N,n))$. Let $(y_1, \ldots, y_N) \in \R^N$ with $0 \leq y_1 < y_2 < \cdots < y_N$.
By expressing an integral as a limit of its Riemann sums, we find from \eqref{eqn:transition} that
\begin{align} \label{eqn:Etransition}
    \Pr(G(n) \in dy \mid G(0) = x) &= \lim_{\ell \to \infty} \Pr(G^{\ell}(n) = y^{\ell} \mid G^{\ell}(0) = x^{\ell}) \ell^N \\
    &= \dt{\partial^{j-i}W_n(y_j-x_i)}_{1\leq i,j \leq N} \, dy.
\end{align}
The last line follows from a direct calculation of the limit using \eqref{eqn:transition}.
The function $W_n$ is the density of the Gamma random variable:
$$ W_n(x) = \frac{x^{n-1}}{(n-1)!}e^{-x} \ind{x \geq 0}.$$
The operator $\partial$ is differentiation $\partial f(x) = f'(x)$ with formal inverse $\partial^{-1} f(x) = \int_{-\infty}^x f(y) dy$.

\subsection{Fixed time law of Exponential lpp with boundary}
Define the following integral kernels on $\R$ for integers $m,n \geq 1$.
\begin{align} \label{eqn:Ekernels}
    Q^m_{E}(u,v) &= \frac{(v-u)^{m-1}}{(m-1)!} e^{u-v} \ind{v \geq u} \\
    S_{E,n,-m}(u,v) &= - \frac{1}{2 \pi \mathbf{i}} \oint_{|w-1|= 1} dw\, e^{(w-1)(v-u)}w^m (1-w)^{-n}\, \ind{v \leq u} \\
    \bar{S}_{E,n,m}(u,v) &= \frac{1}{2 \pi \mathbf{i}} \oint_{|w|= 1} dw\, e^{(w-1)(v-u)} w^{-m} (1-w)^n
\end{align}
Let $B(m)$ be a random walk with $Exp(1)$ step distribution, whose $m$-step transition probability is given by $Q^m_E$ above. Define,
$$ \tau = \inf\, \{m = 0,1,2,3, \ldots: B(m) \leq x_{m+1}\}$$
and
\begin{equation} \label{eqn:Ehypo}
    S^{hypo(x)}_{E,n,m}(u,v) = \E{\bar{S}_{E, n,m-\tau}(B(\tau),v) \ind{\tau < m} \mid B(0) = u}
\end{equation}
Finally, define the integral kernel
\begin{equation} \label{eqn:Ekernel}
    K_E(i,\cdot;j,\cdot) = - Q^{m_j-m_i} \ind{m_i < m_j} + S_{E, n,-m_i} \cdot S^{hypo(x)}_{E, n,m_j}.
\end{equation}
Here $1 \leq m_1,\ldots, m_k$ are distinct integers, $n \geq 1$ and the kernel acts of $L^2(\{1,\ldots,k\} \times \R)$.

\begin{thm} \label{thm:Efixedtime}
    Let $n \geq 1$ and $m_1,\ldots,m_k \geq 1$ be distinct integers. Let $a_1, \ldots, a_k \in \R$. Let $G$ be the corner growth function \eqref{eqn:exppolymer}. Then,
    $$ \Pr(G(m_1,n) \leq a_1, \ldots, G(m_k,n) \leq a_k) = \dt{I - \chi_a K_E \chi_a}_{L^2(\{1,\ldots, k\} \times \R)}$$
    with $\chi_a(i,u) = \ind{u \geq a_i}$.
\end{thm}

\subsection{Proof of Theorem \ref{thm:Efixedtime}}
We take the limit of the formula from Theorem \ref{thm:geomlpp} in the Geometric lpp to Exponential lpp transition.
As stated in \eqref{eqn:geomtoE}, we must perform the limit
$$ \lim_{\ell \to \infty} \Pr(G^{\ell}(m_i,n) \leq a_i; 1 \leq i \leq k)$$
This probability is expressed as a Fredholm determinant by Theorem \ref{thm:geomlpp}.

For a kernel $K$ acting on $\ell^2(\{1,\ldots,k\} \times \Z)$, embed it as a kernel acting on $L^2(\{1, \ldots, k\} \times \R)$ by
\begin{equation}
    K_{\R}(i,u;j,v) = K(i,\lfloor u \rfloor; j, \lfloor v \rfloor)
\end{equation}
One has that
$$ \dt{I+K_{\R}}_{L^2(\{1, \ldots, k\} \times \R)} = \dt{I+K}_{\ell^2(\{1,\ldots,k\} \times \Z)}$$
With this embedding,
$$ \Pr(G^{\ell}(m_i,n) \leq a_i; 1 \leq i \leq k) = \dt{I - \chi^{\ell}K^{\ell} \chi^{\ell}}_{L^2(\{1, \ldots, k\} \times \R)}.$$
where
\begin{equation} \label{eqn:Kell} K^{\ell}(i,u;j,v) = \ell K(i,\lfloor \ell u \rfloor; j, \lfloor \ell v \rfloor).\end{equation}
Here $K$ is the kernel from Theorem \ref{thm:geomlpp} with boundary condition $x^{\ell}_k = \lfloor \ell x_k \rfloor$.
Also, $\chi^{\ell}(i,u) = \ind{\lfloor \ell u \rfloor \geq \ell a_i}$. It is clear that $\chi^{\ell} \to \chi_a$ as $\ell \to \infty$.

We derive the limit of each of the constituent kernels \eqref{eqn:Qgeom}, \eqref{eqn:Sgeom}, \eqref{eqn:Sbargeom} and \eqref{eqn:Shypogeom}
and perform some decay estimates to derive the limit of $K^{\ell}$ as $K_E$.

\subsubsection{Limits of the constituent kernels}
In the formulas \eqref{eqn:Qgeom}, \eqref{eqn:Sgeom}, \eqref{eqn:Sbargeom} and \eqref{eqn:Shypogeom}, choose
$$ \theta = \theta_{\ell} = 1 - \frac{1}{\ell} = q_{\ell} = q.$$
We assume this choice throughout the section.

\begin{lem} \label{lem:QE}
    The kernel $\ell Q^m(\ell u; \ell v)$ converges pointwise in $u,v \in \R$ to $Q^m_{E}(u,v)$ as $\ell \to \infty$.
    Furthermore, it satisfies the decay estimate $\ell | Q^m(\ell u, \ell v)| \leq C_m e^{u-v} |v-u|^{m-1} \ind{v \geq u}$
    for all $u,v \in \R$ and some constant $C_m < \infty$.
\end{lem}

\begin{proof}
Recall $Q(z_1,z_2) = (1-\theta) \theta^{z_2-z_1-1} \ind{z_2 > z_1}$. Thus,
$$ \ell Q(\ell u, \ell v) = \big ( 1 - \frac{1}{\ell} \big)^{\lfloor \ell v \rfloor - \lfloor \ell u \rfloor -1} \ind{\lfloor \ell v \rfloor > \lfloor \ell u \rfloor}.$$
In the limit this goes pointwise to $e^{u-v}\ind{v\geq u} = Q_E(u,v)$. Using the inequality $1+x \leq e^x$, valid for all $x \in \R$, one as
$$ |\ell Q(\ell u, \ell v)| \leq C e^{u-v} \ind{v \geq u}$$
for all $u,v \in \R$. Therefore, one finds that
$| \ell Q^m(\ell u, \ell v)| \leq C_m e^{u-v} (v-u)^{m-1} \ind{v \geq u}$ for all $u,v \in \R$. By the dominated convergence theorem, the pointwise limit, and the decay estimate above, it follows that
$\ell Q^m(\ell u, \ell v) \to Q^m_E(u,v)$ pointwise and also satisfies the decay estimate stated in the lemma.
\end{proof}

\begin{lem} \label{lem:SE}
    The kernel $\ell S_{n,-m}(\ell u, \ell v)$ converges pointwise to $S_{E,n,-m}(u,v)$ as $\ell \to \infty$.
    Furthermore, it satisfies the decay estimate $|\ell S_{n,-m}(\ell u, \ell v)| \leq C_{m,n} e^{v-u} \ind{v \leq u}$ for all $u,v \in \R$ and some constant $C_{m,n}$.
\end{lem}

\begin{proof}
    Looking at the integrand of $S_{n,-m}(u,v)$ from \eqref{eqn:Sgeom}, observe that for large $|w|$ it behaves like $|w|^{z_1-z_2-1-n}$.
    So we can contract the contour to $\infty$ if $z_1-z_2 < n$. Since $m \geq 1$ in the formula of $K$, the poles of the integrand are at $w=0,q$.
    Thus,
    $$ S_{n,-m}(z_1,z_2) = \intz{\gamma_{0,q}} dw\, \frac{\theta^{z_2-z_1}}{w^{z_2-z_1+m+1}} (\frac{1-w}{\alpha})^m \phi(w)^n \; \ind{z_1-z_2 \geq n}.$$
    Here $\gamma_{0,q} = \gamma_0 \cup \gamma_q$ is a union of two small contours around $0$ and $q$.

    Consider the integral over $w \in \gamma_0$. The contour can be contracted to zero if $z_2-z_1+m+1 \leq 0$. So we may assume $z_2-z_1 \geq -m$ in the integral over $\gamma_0$. Then,
    \begin{align*} \ell S_{n,-m}(\ell u, \ell v) &= \frac{\ell}{2 \pi \mathbf{i}} \oint_{\gamma_q} dw\,
    \frac{(1- \frac{1}{\ell})^{\lfloor \ell v \rfloor - \lfloor \ell u \rfloor}}{w^{\lfloor \ell v \rfloor - \lfloor \ell u \rfloor+m+1}} \big (\frac{1-w}{\alpha} \big)^m
    \phi(w)^n\, \ind{\lfloor \ell v \rfloor - \lfloor \ell u \rfloor \leq -n} \\
    & + \frac{\ell}{2 \pi \mathbf{i}} \oint_{\gamma_0} dw \cdot \, \frac{(1- \frac{1}{\ell})^{\lfloor \ell v \rfloor - \lfloor \ell u \rfloor}}{w^{\lfloor \ell v \rfloor - \lfloor \ell u \rfloor+m+1}} \big (\frac{1-w}{\alpha} \big)^m
    \phi(w)^n\, \ind{-m \leq \lfloor \ell v \rfloor - \lfloor \ell u \rfloor \leq -n}.
    \end{align*}
    Since $m,n \geq 1$ in the formula for $K$, the condition $-m \leq \lfloor \ell v \rfloor - \lfloor \ell u \rfloor \leq -n$ is never satisfied for large $\ell$ if $u \neq v$. When $u=v$, the condition $-m \leq 0 \leq -n$ is not met since $m,n \geq 1$.
    So for all $\ell$ sufficiently large, the contribution comes from the integral over $\gamma_q$ only.

    Now change variables $w \mapsto 1 - \frac{w}{\ell}$ in the integral over $\gamma_q$. Let the contour be circular around $q$ with radius $1/\ell$. Since $q = 1-\ell^{-1}$ and $\alpha = (\ell -1)^{-1}$ and $\phi(w) = (1-q)/(w-q)$, one has
    $$\ell S_{n,-m}(\ell u, \ell v) = - \frac{1}{2 \pi \mathbf{i}} \oint_{|w-1|=1} dw\, \frac{e^{u-v} w^m (1-w)^{-n}}{e^{w(u-v)}} \, \ind{v \leq u} \times (1 + O(\ell^{-1})).$$
    Furthermore, using the inequality $|\Re(-\log(1-z)) - \Re(z)| \leq 2 |z|^2$ for $|z| \leq 1/2$, we have the decay estimate
    $$ \ell S_{n,-m}(\ell u, \ell v)| \leq C_{m,n} e^{v-u} \ind{v \leq u}$$
    for all $u,v \in \R$.
\end{proof}

\begin{lem} \label{lem:SbarE}
    The kernel $\ell \bar{S}_{n,m}(\ell u , \ell v)$ converges pointwise to $\bar{S}_{E,n,m}(u,v)$ as $\ell \to \infty$.
    Furthermore, it satisfies the decay estimate $|\ell \bar{S}_{n,m}(\ell u, \ell v)| \leq C_{m,n} e^{u-v}$ for all $u,v \in \R$.
\end{lem}
\begin{proof}
    In the integration contour $|w|=\delta$ from \eqref{eqn:Sbargeom}, choose $\delta = \ell^{-1}$. Make the change of variables $w \mapsto w/\ell$. Then it holds that
    $$\ell \bar{S}_{n,m}(\ell u , \ell v) = \intz{|w|=1}dw\, e^{(w-1)(v-u)} w^{-m} (1-w)^n\, \times (1 + O (\ell^{-1}))$$
    Also, the decay estimate follows as in Lemma \ref{lem:SE} above.
\end{proof}

\begin{lem} \label{lem:SEhypo}
   The kernel $\ell S^{hypo(x^{\ell})}_{n,m}(\ell u, \ell v)$ converges to $S^{hypo(x)}_{E,n,m}(u,v)$ pointwise as $\ell \to \infty$.
   Also, one has the decay estimate $|\ell S^{hypo(x^{\ell})}_{n,m}(\ell u, \ell v)| \leq C_{n,m} e^{-v}$. 
\end{lem}
\begin{proof}
    One has that
    $$\ell S^{hypo(x^{\ell})}_{n,m}(\ell u, \ell v) = \E{\ell \bar{S}_{n,m-\tau_{\ell}}(\ell B^{\ell}(\tau_{\ell}), \ell v) \ind{\tau_{\ell}<m}\mid B^{\ell}(0) = \lfloor \ell u \rfloor /\ell}.$$
    Here $B^{\ell}(m)$ is  a random walk with step distribution $Geom(1-q_{\ell})/\ell$. The hitting time $\tau_{\ell} = \min \, \{m \geq 0: B^{\ell}(m) < \lfloor \ell x_{m+1} \rfloor / \ell\}$.

    There is a coupling of $B^{\ell}$ for all $\ell$ together with a random walk $B$ with $Exp(1)$-step distribution such that $B^{\ell}(m) \to B(m)$ almost surely, simultaneously over all $(\ell, m)$. In this coupling, $\tau_{\ell} \to \tau = \min \, \{ m \geq 0: B(m) \leq x_{m+1}\}$ almost surely.

    Since $| \ell \bar{S}_{n,m}(\ell u, \ell v)| \leq C_{m,n} e^{u-v}$ for all $u,v \in \R$,
    \begin{align*}
        |\ell \bar{S}_{n,m-\tau_{\ell}}(\ell B^{\ell}(\tau_{\ell}), \ell v)| \ind{\tau_{\ell}<m} &\leq C_{m,n} e^{B^{\ell}(\tau_{\ell})-v} \ind{\tau_{\ell} < m} \\
        & \leq C_{m,n} e^{\lfloor \ell x_{\tau_{\ell}+1}\rfloor/\ell -v} \ind{\tau_{\ell}<m} \\
        & \leq C_{n,m} e^{x_{\tau_{\ell}+1} -v + \ell^{-1}} \ind{\tau_{\ell}<m} \\
        & \leq 2C_{n,m} e^{x_{m}-v} \quad (\text{since}\; \tau_{\ell} < m\;\text{and}\; x_{k}\;\text{is increasing}).
    \end{align*}
    By the bounded convergence theorem it follows that for all $u,v \in \R$,
    $$ \ell S^{hypo(x^{\ell})}_{n,m}(\ell u, \ell v) \to S^{hypo}_{E,n,m}(u,v).$$ 

    For the decay estimate, we simply follow the same steps above and use the decay estimate for $\bar{S}_{n,m}$ from Lemma \ref{lem:SbarE} to get that
    \begin{equation*}
        \ell S^{hypo(x^{\ell})}_{n,m}(\ell u, \ell v) \leq C_{m,n} \, e^{x_m-v}  \Pr(\tau^{\ell} < m \mid B^{\ell}(0)=u)  \leq C_{m,n}\, e^{x_m-v}.
    \end{equation*}    
\end{proof}
\subsubsection{Completing the proof} \label{sec:globalproof}
We need to show that the kernel $K^{\ell}$ in \eqref{eqn:Kell} converges pointwise to $K_E$ and satisfies a decay estimate of the form
$$ |K^{\ell}(i,u;j,v)| \leq f_i(u) g_j(v)$$
where $f_i$ are bounded functions and $g_j$ are integrable over $[a_j,\infty)$.

In order to get the decay estimate, we need to append a conjugation factor to the kernel $K^{\ell}$. In integrable probability, it is standard practice to not include conjugation factors in the definition of kernels, but rather introduce them as needed during asymptotics. In fact, even the kernel $K$ from Theorem \ref{thm:geomlpp} should be conjugated in order to be trace class. This is because $Q^m$ has a convolution kernel, and convolution kernels are generally not trace class due to having infinite trace.

Recall that the integers $m_1, \ldots, m_k$ were distinct in the statement of Theorem \ref{thm:geomlpp}.
There is a permutation $\sigma$ of $\{1,\ldots, k\}$ such that $m_{\sigma(1)} < m_{\sigma(2)} < \cdots < m_{\sigma(k)}$.
The inverse $\sigma^{-1}$ then satisfies the property that $m_i < m_j$ if and only if $\sigma^{-1}(i) < \sigma^{-1}(j)$.
Let $\mu = m_{\sigma(k)} - m_{\sigma(1)} = \max_{i,j} \, \{m_j-m_i\}$ be the maximum gap between these integers.

The conjugation factor we append is a multiplication operator:
$$ M f(i,x) = (1+x^2)^{\kappa_{\sigma^{-1}(i)}} f(i,x).$$
The numbers $\kappa_i \geq 0$ need to satisfy $\kappa_j - \kappa_i < -(\mu/2)$ when $i < j$.
For example, we may choose $\kappa_i = \mu (k+1-i)$.
The kernel $K^{\ell}$ is replaced by $M^{-1} K^{\ell} M$ in the determinant, whose entries are thus
\begin{equation} \label{eqn:Kconj}
M^{-1}K^{\ell}M (i,u;j,v) = \frac{(1+v^2)^{\kappa_{\sigma^{-1}(j)}}}{(1+u^2)^{\kappa_{\sigma^{-1}(i)}}} K^{\ell}(i,u;j,v).
\end{equation}

\begin{lem} \label{lem:Qdecay}
    Suppose $m_i < m_j$, so that $\sigma^{-1}(i) = r < \sigma^{-1}(j) = s$. Let $m = m_j-m_i$.
    There are functions $f_i$ and $g_j$ such that $f_i$ is bounded over $[a_i,\infty)$, $g_j$ is integrable over $[a_j, \infty)$, and
    $$\frac{(1+v^2)^{\kappa_s}}{(1+u^2)^{\kappa_r}} |(\ell Q^m(\ell u, \ell v)| \leq f_i(u) g_j(v).$$
\end{lem}
\begin{proof}
    From Lemma \ref{lem:QE} we see that $| \ell Q^m(\ell u, \ell v)| \leq C_{m} e^{u-v} (v-u)^{m-1} \ind{v \geq u}$.
    We have $2 (\kappa_s-\kappa_r) < -\mu \leq -m$. Thus, there is a number $\kappa$ satisfying $2\kappa_s + m < \kappa \leq 2\kappa_r$.

    Since $u \geq a_i$, $|v-u|^{m-1} \leq C_{a_i} |v|^{m-1}$.
    Write $e^{u-v} = e^{u-v - \kappa \log|v| + \kappa \log|v|} = |v|^{-\kappa} e^{u-v+\kappa \log |v|}$.
    Set $g_j(v) = (1+v^2)^{\kappa_s} |v|^{m-1 - \kappa}$. Now  $2\kappa_s +m-1 - \kappa < -1$, and so $g_j(v)$ is integrable over $v \in [a_j,\infty)$.
    Set $f_i(u) = (1+u^2)^{-\kappa_r} \sup_{v \geq a_j} e^{u-v +\kappa \log |v|} \ind{v \geq u}$. The supremum is obtained at $v = u + O(\kappa)$.
    Therefore, $|f_i(u)| \leq C_{\kappa} |u|^{-2\kappa_r +\kappa}$, and this is bounded due to $\kappa \leq 2\kappa_r$.
\end{proof}

\begin{lem} \label{lem:SSbar}
    The kernel $\ell S_{n,-m_i} \cdot S^{hypo}_{n,m_j}(\ell u, \ell v)$ converges pointwise to $S_{E,n,-m_i} \cdot S^{hypo(x)}_{E,n,m_j}(u,v)$.
    Furthermore, one has the estimate: $|\ell S_{n,-m_i} \cdot S^{hypo}_{n,m_j}(\ell u, \ell v)| \leq C e^{-v}$ for a constant $C$ free of $\ell$, $u$ and $v$.
\end{lem}

\begin{proof}
    One has
    \begin{equation} \label{eqn:SSbar1}
    \ell S_{n,-m_i} \cdot S^{hypo}_{n,m_j}(\ell u, \ell v) = \int_{\R} dz\, \ell S_{n,-m_i}(\ell u, \ell z) \cdot \ell S^{hypo(x^{\ell}}_{n,m_j}(\ell z, \ell v).
    \end{equation}
    From Lemmas \ref{lem:SE} and \ref{lem:SEhypo}, it follows that the integrand converges pointwise to $S_{E, n,-m_i}(u,z) \cdot S^{hypo(x)}_{E,n,m_j}(z,v)$.
    Furthermore, one has the decay estimate,
    $$| \ell S_{n,-m_i}(\ell u, \ell z) \cdot \ell S^{hypo(x^{\ell}}_{n,m_j}(\ell z, \ell v)| \leq C e^{z-u + x_{m_j}-v} \ind{z \leq u}.$$
    The right side is integrable over $z \in \R$. The dominated convergence theorem implies
    $$\ell S_{n,-m_i} \cdot S^{hypo}_{n,m_j}(\ell u, \ell v) \to S_{E,n,-m_i} \cdot S^{hypo(x)}_{E,n,m_j}(u,v)$$
    for every $u,v \in \R$.

    Finally, the decay estimate above also shows
    \begin{align*}
        |\ell S_{n,-m_i} \cdot S^{hypo}_{n,m_j}(\ell u, \ell v)| & \leq \int_{\R}dz\, C e^{z-u + x_{m_j}-v} \ind{z \leq u} \\
        & \leq C' e^{m_j-v}.
    \end{align*}
\end{proof}

Lemmas \ref{lem:QE} and \ref{lem:Qdecay} imply that the $\ell Q^{m_j-m_i}(\ell u, \ell v) \ind{m_i < m_j}$ converges to $Q^{m_j-m_i}_E \ind{m_i < m_j}$
pointwise and satisfies a sufficient global decay estimate after conjugation. Lemma \ref{lem:SSbar} says $\ell S_{n,-m_i} \cdot S^{hypo}_{n,m_j}(\ell u, \ell v)$ converges pointwise to $S_{E,n,-m_i} \cdot S^{hypo(x)}_{E,n,m_j}(u,v)$. After conjugation, its global decay is of order
$$C_{m_i,m_j,n}\, \frac{(1+v^2)^{a_{\sigma^{-1}(j)}}}{(1+u^2)^{a_{\sigma^{-1}(i)}}} e^{-v}$$
which is integrable of $v \in [a_j,\infty)$ and bounded in $u \in [a_i,\infty)$ as the $a_i \geq 0$. This completes the proof of Theorem \ref{thm:Efixedtime}.

\subsection{Continuous time tasep} \label{sec:ctasep}
Continuous time tasep in an interacting particle system consisting of particles on $\Z$ that move in continuous time $t \geq 0$
There is some initial configuration of particles at time $t = 0$. Each particle then attempts to jump one unit rightward at rate 1.
The jump is successful if there is no particle already occupying the jump site (the exclusion rule).
See \cite{Liggett} for the fact that this process is probabilistically well defined.

Let us focus on the case when there is a rightmost particle, which is particle number $0$.
Label the particles from right to left by integers $n = 0,1,2,3, \cdots$. Let $Y_n(t)$ be the location
of particle number $n$ at time $t$. Let
$$ G_Y(m,n) = \text{time when particle number}\; n\; \text{makes jump number}\; m.$$
for $n \geq 0$ and $m \geq 1$. The process $G_Y$ and the initial condition $Y_n(0)$
determine the evolution of tasep. The function $G_Y$ satisfies the recursion
\begin{equation} \label{eqn:etasep}
G_Y(m,n) = \max \{G_Y(m-1,n), G_Y(m+1-\mathrm{gap}_n, n-1) \} + \omega_{m,n}
\end{equation}
where $\omega_{m,n}$ are independent $Exp(1)$ random variables and $\rm{gap}_n = Y_{n-1}(0) - Y_n(0) \geq 1$.

Consider the step initial condition with particles initially at sites $-n$ for $n \geq 0$. Particle number $n$ is the one initially at site $-n$.
Let $X^*_n(t)$ be the position of particle number $n$ at time $t$. This process is the space-noise dual of the Exponential lpp model.
The boundary condition $0 < x_1 < x_2 < \cdots$ translates to particle $X^*_0$ having a deterministic trajectory:
it jumps one unit rightward at the times $x_1, x_2$, etc. The other particles move randomly at rate 1. Theorem \ref{thm:Efixedtime} has the following corollary, by way of property (7) of the space-noise dual.

\begin{cor} \label{cor:ctasep}
    Consider continuous time tasep from the step initial condition $X^*_n(0) = -n$ for $n \geq 0$.
    Let $X^*_n(t)$ be the location of particle number $n$ at time $t$. Let $X^*_0$ jump one unit rightward at the times $x_k \in \R$ with $0 < x_1 < x_2 < x_3 < \cdots$. The other particles move randomly at rate 1 subject to the exclusion rule. Then for $n \geq 1$, times $t_1, t_2, \ldots, t_k > 0$, and distinct integers $a_1, \ldots, a_k \geq -n$,
    \begin{equation} \label{eqn:ctasepformula2}
    \Pr(X^*_n(t_i) > a_i; 1 \leq i \leq k) = \dt{I - \chi_t J_E \chi_t}_{L^2(\{1,2,\ldots,k\}\times \R)}
    \end{equation}
    where $\chi_t(i,z) = \ind{z \geq t_i}$ and the kernel $J_E$ is the kernel $K_E$ from Theorem \ref{thm:Efixedtime} with choice of parameters $m_i = a_i+n+1$ and $n$ the same.
\end{cor}

\section{Brownian lpp and Brownian tasep} \label{sec:brown}
Brownian lpp operates on the domain $D_s = [0,\infty)$, $D_t = \{1,2,3,\ldots\}$ and $D_\eta = \R$.
The noise is given by a collection of independent, standard Brownian motions $x \mapsto \eta(x,i)$, for $i \geq 1$.
Introduce a boundary condition with a continuous function $b : [0,\infty) \to \R$ with $b(0) = 0$.
Extend the temporal domain to $\{0,1,2,3,\ldots\}$ and set $\eta(x,0) = b(x)$. Consider the corner growth function $G$ from $(0,0)$:
\begin{equation} \label{eqn:GB}
    G(x,n) = L_B(0,0;x,n) = \max_{0 \leq x_0 \leq \cdots \leq x_n = x} b(x_0) + \sum_{i=1}^n \eta(x_i,i) - \eta(x_{i-1},i).
\end{equation}

\subsection{Exponential to Brownian transition}
Consider the Exponential lpp model \eqref{eqn:Enoise}. It can be used to transition into the Brownian lpp model, as originally observed by Glynn and Whitt \cite{GW}. Extend the noise $\eta(x,i)$ in Exponential lpp to $x \in [0,\infty)$ by linear interpolation ($i \geq 0$ remains discrete). 
The boundary condition is $\eta(x,0) = x + \sqrt{N}b(x/N)$, where $N \geq 1$ be a large integer scaling parameter. Define
\begin{align}
    \eta_N(x,i) &= \frac{\eta(Nx,i)-Nx}{\sqrt{N}} \quad i \geq 1,\\
    \eta_N(x,0) &= \frac{\eta(Nx,0) -Nx}{\sqrt{N}} = b(x).
\end{align}
Denote by $G_E$ the corner growth function of the Exponential lpp model \eqref{eqn:Elpprecur} with boundary condition $x_k = k + \sqrt{N} b(k/N)$.
The scaled noise $\eta_N$ is the noise of an lpp model on the domain $D_s = [0,\infty)$, $D_t = \{0,1,2,\ldots\}$ and $D_{\eta} = \R$.
Its corner growth function $G_N$ from $(0,0)$ satisfies
\begin{equation} \label{eqn:GN}
    G_N(x,n) = \max_{0 \leq x_0 \leq \cdots \leq x_n = x} b(x_0) + \sum_{i=1}^n \eta_N(x_i,i) - \eta_N(x_{i-1},i) = \frac{G_E(\lceil Nx \rceil,n)-Nx-n}{\sqrt{N}} + O(n/N).
\end{equation}

By Donsker's theorem, there exists a collection $B_i(x)$ of independent, standard Brownian motions coupled with the noise $\eta_N$ for every $N$ such that
$$ \eta_N(x,i) \to B_i(x) \quad \text{almost surely}$$
uniformly over compact subsets of $(x,i) \in [0,\infty) \times \mathbb{N}$. In fact, since the increments of $\eta(x,i)$ in Exponential lpp are Exponential random variables, one has a quantitative rate of convergence via the KMT approximation \cite{Chat}. For every real $T > 0$ and integer $n \geq 1$, there is a random constant $C_{T,n}$
with $\E{C_{T,n}} < \infty$ such that
\begin{equation} \label{eqn:KMT} \sup_{x \in [0,T], 1 \leq i \leq n} |\eta_N(x,i) - B_i(x)| \leq C_{T,n} \frac{\log N}{\sqrt{N}}. \end{equation}

From this convergence of the noise it follows that the lpp model with noise $\eta_N$ converges to the Brownian lpp model. In particular, for the corner growth function $G(x,n)$ in \eqref{eqn:GB}, one has
$$ G(x,n) = \lim_{N \to \infty} G_N(x,n)$$
almost surely (here the Brownian lpp is driven by the Brownian motions in the aforementioned coupling), with the convergence being uniformly over compacts.
Consequently, for $x_1, \ldots, x_k > 0$, $n \geq 1$ and $a_1, \ldots, a_k \in \R$, due to \eqref{eqn:GN},
\begin{equation} \label{eqn:EBlimit}
    \Pr(G(x_i,n) \leq a_i, 1 \leq i \leq k) = \lim_{N \to \infty} \, \Pr(G_E(\lceil Nx_i \rceil, n) \leq Nx_i+n+\sqrt{N} a_i, 1 \leq i \leq k)
\end{equation}

\subsection{Fixed time law of Brownian lpp with boundary}
Define the following integral kernels acting on $\R$.

For $x > 0$, define the heat kernel
$$e^{\frac{x}{2} \partial^2}(u,v) = \frac{1}{\sqrt{2\pi x}} e^{-\frac{(u-v)^2}{2x}}.$$
For $x > 0$ and $n \geq 1$, define
\begin{align}
S_{B,n,-x}(u,v) &= \frac{1}{2\pi \mathbf{i}} \oint_{|w|=1} dw\, e^{-\frac{x}{2}w^2 + (u-v)w} w^{-n},\\
\bar{S}_{B,n,x}(u,v) &= \frac{1}{2\pi \mathbf{i}} \oint_{\Re(w)=0} dw\, e^{\frac{x}{2}w^2 + (u-v)w} w^{n}.
\end{align}
The contour $\{\Re(z)=0\} = \{\mathbf{i}t; t \in \R\}$ is oriented upward.

These kernels can be expressed in terms of Hermite polynomials.
Recall the Hermite polynomials
\begin{equation} \label{eqn:Hermite}
 H_n(x) = (-1)^n e^{x^2/2} \partial^n e^{-x^2/2}, \quad e^{xz - \frac{1}{2}z^2} = \sum_{n=0}^{\infty} H_n(x) \frac{z^n}{n!}.
\end{equation}
Then,
\begin{align} \label{eqn:SH}
    S_{B,n,-x}(u,v) &= \frac{x^{(n-1)/2}}{(n-1)!} H_{n-1}\big(\frac{u-v}{\sqrt{x}}\big), \\
    \bar{S}_{B,n,x}(u,v) &= (-1)^n x^{-n/2} e^{\frac{x}{2} \partial^2}(u,v) \, H_n \big ( \frac{u-v}{\sqrt{x}} \big).
\end{align}
Define the kernel
\begin{equation} \label{eqn:SBhypo} S^{hypo(b)}_{B, n,x}(u,v) = \E{\bar{S}_{B, n,x-\tau}(B(\tau),v) \ind{\tau \leq x} \mid B(0)=u}\end{equation}
where $B(y)$ is a standard Brownian motion started from $B(0)=u$ and
$$\tau = \inf \{y \geq 0: B(y) \leq b(y)\}.$$

The kernel $S^{hypo(b)}_{B,n,x}$ can be interpreted in terms of a hitting probability. We have that $\bar{S}_{B,n,x}(u,v) = \partial^n e^{\frac{x}{2}\partial^2}(u,v)$.
Since $e^{\frac{x}{2}\partial^2}(u,v) = e^{(u-v)^2/2x}/\sqrt{2 \pi x}$,
$$ \partial^n e^{\frac{x}{2}\partial^2}(u,v) = \partial_u^n e^{(u-v)^2/2x}/\sqrt{2 \pi x} = (-1)^n \partial^n_v e^{(u-v)^2/2x}/\sqrt{2 \pi x}.$$
As a result, $S^{hypo(b)}_{B,n,x}(u,v) = (-\partial_v)^n \E{e^{(\frac{x-\tau)}{2}\partial^2}(B(\tau),v) \ind{\tau \leq x}\mid B(0)=u}$.
By the strong Markov property of Brownian motion,
$$\E{e^{(\frac{x-\tau)}{2}\partial^2}(B(\tau),v) \ind{\tau \leq x}\mid B(0)=u} = \Pr(\tau \leq x, B(x) \in dv \mid B(0)=u).$$
Therefore,
\begin{equation} \label{eng:SBhypoprob}S^{hypo(b)}_{B,n,x}(u,v) = (-\partial_v)^n \Pr(\tau \leq x, B(x) \in dv \mid B(0)=u).\end{equation}

Finally, for $n \geq 1$ and distinct $x_1, \ldots, x_k > 0$, define the kernel $K_B$ according to
\begin{equation}\label{eqn:KBrownian}
K_B(i, \cdot; j, \cdot) = - e^{\frac{(x_j-x_i)}{2}\partial^2} \mathbf{1}_{x_i < x_j} + S_{B,n, -x_i}S^{hypo(b)}_{B,n,x_j}.
\end{equation}
The kernel acts on the space $L^2(\{1,\ldots, k\} \times \mathbb{R})$.

\begin{thm} \label{thm:Blpp}
Let $n \geq 1$ be an integer, $x_1, \ldots, x_k > 0$ and $a_1, \ldots, a_k$ be real numbers such that the $x_i$ are distinct.
Assume $b : [0,\infty] \to \mathbb{R}$ is continuous and $b(0) = 0$. For the growth function $G$ in \eqref{eqn:GB},
$$Pr(G(x_i,n) \leq a_i; 1 \leq i \leq k) = \det(I - \chi_a K_B \chi_a)_{L^2(\{1,\ldots, k\} \times \mathbb{R})}$$
with $\chi_a(i,z) = \ind{z \geq a_i}$.
\end{thm}

\begin{exm} \label{exm:flat}
    Consider the boundary $b \equiv 0$. In order to compute the kernel $K_B$, we need to calculate $P_x(u,v) = \Pr(\tau \leq x, W(x) \in dv \mid W(0) = u)$
    where $x > 0$ and $W$ is a standard Brownian motion. The event $\{\tau \leq x\}$ occurs automatically if $u \leq 0$ or $u > 0$ and $v \leq 0$.
    For $u,v > 0$, the reflection principle gives $P_x(u,v) = e^{\frac{x}{2}\partial^2}(-u,v)$. Therefore,
    \begin{align*}
    S^{hypo(b)}_{B,n,x}(u,v) &= (-\partial_v)^n P_x(u,v) = x^{-n/2} H_n(\frac{v-u}{\sqrt{x}}) e^{-\frac{(u-v)^2}{2x}} \left( \mathbf{1}_{u \leq 0} + \mathbf{1}_{u > 0, v \leq 0}\right)\\
    &+ x^{-n/2} H_n(\frac{v+u}{\sqrt{x}}) e^{-\frac{(u+v)^2}{2x}} \mathbf{1}_{u>0,v>0}.
    \end{align*}

    Suppose $v \leq 0$. Then,
    \begin{align*}
       S_{B,n,-x_1} \cdot S^{hypo(b)}_{B, n,x_2}(u,v) &= \frac{x_1^{(n-1)/2} x_2^{-n/2}}{(n-1)! \sqrt{2\pi  x_2}} \int_{-\infty}^{\infty} dz\,H_{n-1}(\frac{u-z}{\sqrt{x_1}}) H_n(\frac{v-z}{\sqrt{x_2}}) e^{\frac{(v-z)^2}{2x_2}} \\
       &= \frac{x_1^{(n-1)/2} x_2^{-n/2}}{(n-1)! \sqrt{2\pi  x_2}} \int_{-\infty}^{\infty} dz\, H_{n-1}\big(\frac{u-v+ \sqrt{x_2}z}{\sqrt{x_1}}\big) \, H_n(z) e^{-z^2/2}
    \end{align*}
    The integral above is zero because $H_{n-1}(\cdot)$ is a polynomial in $z$ or degree at most $n-1$ and $H_n$ is orthogonal to such polynomials under the density $e^{-z^2/2}$.

    Suppose $v > 0$. A calculation similar to the above gives:
    $$ S_{B,n,-x_1} \cdot S^{hypo(b)}_{B, n,x_2}(u,v)  = \frac{x_1^{(n-1)/2} x_2^{-n/2}}{(n-1)! \sqrt{2\pi  x_2}}
    \int_0^{\infty} dz\, \Big ( H_{n-1}\big(\frac{u+z}{\sqrt{x_1}}\big) + H_{n-1}\big(\frac{u-z}{\sqrt{x_1}}\big)\Big) H_n(\frac{z+v}{\sqrt{x_2}}) e^{-\frac{(z+v)^2}{2x_2}}.$$
    As such, we find that
    \begin{align*}
        K_B(i,u;j,v) &= - e^{-\frac{(u-v)^2}{2(x_j-x_i)}} \mathbf{1}_{x_i < x_j} + \mathbf{1}_{ v > 0} \, \frac{x_i^{(n-1)/2} x_j^{-n/2}}{(n-1)! \sqrt{2\pi  x_j}} \times \\
        &\int_0^{\infty} dz\, \Big ( H_{n-1}\big(\frac{u+z}{\sqrt{x_i}}\big) + H_{n-1}\big(\frac{u-z}{\sqrt{x_i}}\big)\Big) H_n\big(\frac{z+v}{\sqrt{x_j}}\big) e^{-\frac{(z+v)^2}{2x_j}}.
    \end{align*}
\end{exm}

\subsection{Proof of Theorem \ref{thm:Blpp}}
We use Theorem \ref{thm:Efixedtime} to take the the limit on the right side of \eqref{eqn:EBlimit} and obtain a Fredholm determinant.
Consider the Exponential lpp model with boundary conditions $x_k = k + \sqrt{N}b(k/N)$. Then,
$$\Pr(G_E(\lceil Nx_i \rceil, n) \leq Nx_i+n+\sqrt{N} a_i, 1 \leq i \leq k) = \dt{I - \chi^N K^N \chi^N}_{L^2(\{1,\ldots,k\} \times \R)}$$
where $\chi^N(i,u) = \ind{u \geq a_i - (n/\sqrt{N})}$ and
$$K^N(i,u;j,v) = \sqrt{N} K_E(i,Nx_i +\sqrt{N}u ;j, Nx_j + \sqrt{N}v).$$
It is clear that $\chi^N \to \chi_a$ as $N \to \infty$.

As in the proof of Theorem \ref{thm:Efixedtime}, we derive the limit of each of the constituent kernels \eqref{eqn:Ekernels} and  \eqref{eqn:Ehypo}, and then perform some decay estimates to derive the limit of $K^{N}$ as $K_B$.

\subsubsection{Limits of constituent kernels}
\begin{lem} \label{lem:QBlimit}
    Suppose $x_1 < x_2$. Let $n_1 = \lceil Nx_1 \rceil$, $n_2 = \lceil Nx_2 \rceil$, $u_1 = Nx_1 + \sqrt{N}u$ and $u_2 = Nx_2 + \sqrt{N} v$. Then,
    $$ \sqrt{N} Q_E^{n_2-n_1}(u_1,u_2) \to e^{\frac{(x_2-x_1)}{2}\partial^2}(u,v)$$
    pointwise in $u,v \in \R$. Furthermore, one has the decay estimate: $ |\sqrt{N} Q_E^{n_2-n_1}(u_1,u_2)| \leq C_{x_1,x_2} e^{\frac{(u-v)^2}{2(x_2-x_1)}}$ for all $u,v \in \R$.
\end{lem}

\begin{proof}
    We may use Stirling's approximation of factorials: $n! = \sqrt{2 \pi n} (n/e)^n \times (1 + O(1/n))$. The error term is uniformly bounded in $n$, so there is a constant $C$ such that $n! \leq C \sqrt{2\pi n} (n/e)^n$.

    Recall $Q^m(u_1,u_2) = e^{u_1-u_2} \frac{(u_2-u_1)^{m-1}}{(m-1)!} \ind{u_2 \geq u_1}$. Using Stirling's approximation for $(m-1)!$, $m = n_2-n_1$, and some bookkeeping, one finds $$ \sqrt{N} Q^m(u_1,u_2) = (\sqrt{2 \pi (x_2-x_1)})^{-1} e^{(v-u)^2/2(x_2-x_1)} \times (1 + O(1/N)).$$ Also, the multiplicative error term is uniformly bounded by some constant $C$ for all $u,v \in \R$.
\end{proof}

\begin{lem} \label{lem:SBlimit}
  Suppose $x > 0$, $m = \lceil Nx \rceil$ and $n \geq 1$. Let $u_1 = Nx + \sqrt{N} u$ and $u_2 = \sqrt{N} v$. Then,
  $$ \sqrt{N} S_{E,n,-m}(u_1,u_2) = S_{B, n, -x}(u,v) \times \ind{v-u \leq \sqrt{N} x}\,  N^{n/2} (1 + O(x N^{-1/2})).$$
  The constant in the error term $O(xN^{-1/2})$ is uniformly bounded in $N \geq 1$, $u,v \in \R$ and $x \geq 0$.
\end{lem}

\begin{proof}
    From the definition, we have
    \begin{align*}
        \sqrt{N} S_{E,n,-m}(u_1,u_2) &= -\frac{\sqrt{N}}{2 \pi \mathbf{i}} \oint_{|w-1|=1} dw\, e^{(w-1)(\sqrt{N}(v-u)- Nx)} w^{\lceil Nx \rceil} (1-w)^{-n} \; \ind{v-u \leq \sqrt{N}x}\\
        &=-\frac{\sqrt{N}}{2 \pi \mathbf{i}} \oint_{|w-1|=\delta > 0} dw\, e^{(w-1)(\sqrt{N}(v-u)- Nx)} w^{\lceil Nx \rceil} (1-w)^{-n} \; \ind{v-u \leq \sqrt{N}x}.
    \end{align*}
    Change variables $w = 1 - \frac{z}{N}$ and $\delta = \delta'/\sqrt{N}$ for $\delta' > 0$. Then the above becomes
    \begin{equation*}
       \sqrt{N} S_{E,n,-m}(u_1,u_2) =  \frac{1}{2 \pi \mathbf{i}} \oint_{|z|=\delta' > 0} dz\, e^{z(u-v) + \sqrt{N}z x + \lceil Nx \rceil \log(1 - z/N)} (z/\sqrt{N})^{-n} \; \ind{v-u \leq \sqrt{N}x}.
    \end{equation*}
    If $|z| = \delta' \leq 1/2$ then Taylor expansion gives $\log(1 - \frac{z}{N}) = - \frac{z}{\sqrt{N}} - \frac{z^2}{2N} + O(N^{-3/2})$ where the big O term has a absolute constant (free of all parameters). Since $\lceil Nx \rceil - Nx \in [0,1]$, we find that
    \begin{equation*}
        \sqrt{N} S_{E,n,-m}(u_1,u_2) = \frac{1}{2 \pi \mathbf{i}} \oint_{|z|=\delta'} dz\, e^{- \frac{x}{2}z^2 + (u-v)z} z^{-n} \times N^{n/2} (1 + O(xN^{-1/2})) \ind{v-u \leq \sqrt{N}x}.
    \end{equation*}
    We recognize the integral on the right side as $S_{B,n,-x}(u,v)$.
\end{proof}

\begin{lem} \label{lem:barSBlimit}
    Suppose $x > 0$, $m = \lceil Nx \rceil$ and $n \geq 1$. Let $u_2 = Nx + \sqrt{N} v$ and $u_1 = \sqrt{N} u$. Then,
  $$ \sqrt{N} \bar{S}_{E,n,m}(u_1,u_2) = \bar{S}_{B, n, x}(u,v)\,  N^{-n/2} (1 + O(x N^{-1/2})).$$
  The constant in the error term $O(xN^{-1/2})$ is uniformly bounded in $N \geq 1$, $u,v \in \R$ and $x \geq 0$.
\end{lem}
\begin{proof}
    By definition of $\bar{S}_E$,
    \begin{equation*}
        \sqrt{N} \bar{S}_{E,n,m}(u_1,u_2) = \frac{\sqrt{N}}{2 \pi \mathbf{i}} \oint_{|w|=\delta} dw\, e^{(w-1)(\sqrt{N}(v-u)+ Nx)} w^{-\lceil Nx \rceil} (1-w)^{n}.
    \end{equation*}
    
    We choose $\delta = 1$ and parametrize $w = w(\theta) = e^{\mathbf{i} \theta / \sqrt{N}}$ for $|\theta| \leq \pi \sqrt{N}$. We change variables
    $w = 1 - \frac{z}{\sqrt{N}}$.  Then $z(\theta) = (1 - e^{\mathbf{i} \theta/\sqrt{N}}) \sqrt{N}$.
    Let $\gamma_N$ denote the closed contour sketched out by $z(\theta)$. We have $|z(\theta)| \leq \sqrt{N} \min\{2, |\theta|/\sqrt{N}\}$ and,
    locally, $z(\theta) = -\mathbf{i} \theta$ if $|\theta| = O(1)$ in $N$.

    In the new variable we find that
    $$\sqrt{N} \bar{S}_{E,n,m}(u_1,u_2) = - \frac{1}{2 \pi \mathbf{i}} \oint_{\gamma_N} dz\, e^{F_{N,x}(z) + z(u-v)} z^n \times N^{-n/2}$$
    with $F_{N,x}(z) = -\sqrt{N} x z - Nx \log(1 - z/\sqrt{N})$.

    Locally, if $|\theta| \leq L$, then by Taylor expansion $F_{N,x}(z(\theta)) = \frac{x}{2} (z(\theta)^2 + O(x L / \sqrt{N})$ and $z(\theta) = -\mathbf{i} \theta + O(L/\sqrt{N})$. The big O constant is free of all parameters.
    The contour $\gamma_N$ approximates the vertical line $\{ \Re(z) = 0\} = \{\mathbf{i} \theta: \theta \in \R\}$, oriented downwards.
    The integrand converges pointwise in $\theta$ to $e^{\frac{x}{2} z^2 + z(u-v)} z^n$ with $z = z(\theta) = -\mathbf{i} \theta$.

    Globally, since $1- \frac{z(\theta)}{\sqrt{N}} = e^{ \mathbf{i} \theta /\sqrt{N}}$, we have $\Re(F_{N,x}(z(\theta))) = - N x \Re(z(\theta)) -Nx \log(1) = -Nx (1- \cos(\theta/\sqrt{N}))$. For $|\theta| \leq \pi \sqrt{N}$, $1 - \cos(\theta/\sqrt{N}) \geq \theta^2 / 10 N$. Therefore, $\Re(F_{N,x}(z(\theta)) \leq - \theta^2 x/10$ since $x \geq 0$. Also, $\Re(z(\theta)(u-v)) \leq |\theta| |u-v|$ because $|z(\theta)| \leq |\theta|$. Similarly, $|z(\theta)|^n \leq |\theta|^n$.
    Therefore, the integrand is bounded above by $e^{-\theta^2 x / 10 + |\theta||u-v| + n \log |\theta|}$.

    From the dominated convergence theorem, and re-orienting the limit contour upwards, we have that
    $$ \sqrt{N} \bar{S}_{E,n,m}(u_1,u_2) = \bar{S}_{B, n,x}(u,v) \, N^{-n/2} \times (1 + O(x/\sqrt{N})).$$
\end{proof}

\begin{lem} \label{lem:SBhypolimit}
    Suppose $x > 0$, $m = \lceil Nx \rceil$ and $n \geq 1$. Let $u_2 = Nx + \sqrt{N} v$ and $u_1 = \sqrt{N} u$. Let $y^N_k = k + \sqrt{N}b(k/N)$.
    Then,
    $$ \sqrt{N} S^{hypo(y^N)}_{E,n,m}(u_1,u_2) = S^{hypo(b)}_{B,n,x}(u,v) \times N^{-n/2} (1 + O_{x,n}(1/\sqrt{N})).$$
\end{lem}

\begin{proof}
    We have
    $$\sqrt{N} S^{hypo(y^N)}_{E,n,m}(u_1,u_2) = \E{\sqrt{N} \bar{S}_{E,n,\lceil Nx \rceil - \tau_N}(B_n(\tau_N), Nx + \sqrt{N}v)\ind{\tau_N < \lceil Nx \rceil}\mid B_N(0) = \sqrt{N}u}.$$
    Here $B_N(m)$ is a random walk with $Exp(1)$ step distribution and $\tau_N = \inf \{m\geq 0: B_n(m) \leq y^N_{m+1}\}$.

    Define the scaled walk $\bar{B}_N(y) = \frac{B_N(\lceil Ny \rceil)-Ny}{\sqrt{N}}$ and $\bar{\tau}_N = \tau_N/N$.
    Using the fact that $\bar{S}_{E,n,m}(u+a,v) = \bar{S}_{E,n,m}(u,v-a)$, and Lemma \ref{lem:barSBlimit}, we see that
    $$\sqrt{N} S^{hypo(y^N)}_{E,n,m}(u_1,u_2) = \E{\bar{S}_{B,n,x-\bar{\tau}_N}(\bar{B}_N(\bar{\tau}_N),v) \ind{\bar{\tau}_N \leq x} \mid \bar{B}_N(0)=u} N^{-n/2} (1+ O(x/\sqrt{N})).$$

    By Donsker's theorem and Shorokhod representation theorem, we may find a coupling of the scaled walks $\bar{B}_N$ for every $N$ together with a Brownian motion $B$ such that $\bar{B}_N \to B$ uniformly on compacts. In fact, the KMT approximation gives a coupling satisfying
    $$ \sup_{y \in [0,x]} |\bar{B}_N(y) - B(y)| \leq C_x (\log N)/\sqrt{N}$$
    where $\E{C_x} < \infty$. From this and the smoothness of $\bar{S}_{B, n,x}(u,v)$ in all parameters, we can deduce that
    $$ \big |\E{\bar{S}_{B,n,x-\bar{\tau}_N}(\bar{B}_N(\bar{\tau}_N,v)) \ind{\bar{\tau}_N \leq x} \mid \bar{B}_N(0)=u} - \E{\bar{S}_{B,n,x-\bar{\tau}}(B(\bar{\tau},v)) \ind{\bar{\tau} \leq x} \mid B(0)=u} \big| \leq C \frac{\log N}{N^{1/2}}$$
    where $\bar{\tau} = \inf \{y \geq 0: B(y) \leq b(y)\}$ is the almost sure limit of $\bar{\tau}_N$. The lemma follows.
\end{proof}

\subsubsection{Completing the proof}
Lemmas \ref{lem:QBlimit}, \ref{lem:SBlimit}, and \ref{lem:SBhypolimit} establish pointwise convergence of the constituent kernels of $K_B$.
To complete the proof, we have to establish decay estimates so that the Fredholm series expansion converges absolutely.
That the errors provided by the aforementioned lemmas are uniformly bounded in all parameters.
Thus, it is enough to show decay estimates for the kernels $Q_{B}$, $S_{B,n,-x}$ and $S^{hypo(b)}_{B,n,x}$.
This will also confirm that $K_B$ has an absolutely convergent Fredholm series.

In order to get decay estimates, we have to conjugate the kernel $K_B$. The conjugation factor we need for the $(i,j)$ block is
\begin{equation} \label{eqn:conjugation}
    \frac{(1+v^2)^{k+n-\sigma^{-1}(j)}}{(1+u^2)^{k+n-\sigma^{-1}(i)}}
\end{equation}
where $\sigma$ is the permutation for which $x_{\sigma(1)} < x_{\sigma(2)}  < \cdots < x_{\sigma(k)}$.

With this conjugation factor, arguing as in the proof of Lemma \ref{lem:Qdecay}, it follows that
$$\frac{(1+v^2)^{k+n-\sigma^{-1}(j)}}{(1+u^2)^{k+n-\sigma^{-1}(i)}} \, e^{\frac{x_j-x_i}{2} \partial^2}(u,v) \ind{x_i < x_j} \leq f_i(u) g_j(v)$$
where $f_i$ is bounded and $g_j$ is integrable over $\R$.

\begin{lem} \label{lem:Hermitedecay}
    Suppose $n \geq 1$ and $x \in [0,T]$. There is a constant $C_{n,T}$ such that
    $$ |S_{B.n,-x}(u,v)| \leq C_{n,T} ( |u-v|^{n-1} + 1) \quad \text{and} \quad |\bar{S}_{B,n,x}(u,v)| \leq C_{n,T} (|u-v|^n+1) e^{\frac{x}{2}\partial^2}(u,v).$$
\end{lem}
\begin{proof}
    The bounds follow from \eqref{eqn:SH}, which represents these kernels in terms of Hermite polynomials.
\end{proof}

\begin{lem} \label{lem:tauprob}
    Let $B$ be a standard Brownian motion, $x > 0$, and $\tau = \inf\{y\geq 0: B(y) \leq b(y)\}$. Let $\mu = \max_{y \in [0,x]} b(y)$. If $z > \mu$ then
    $\Pr(\tau \leq x \mid B(0)=z) \leq 2 e^{- (z-\mu)^2/2x}$.
\end{lem}

\begin{proof}
    In order for $B$ to hit $b$ by time $x$, it must go below level $\mu$ since $z > \mu$. The latter event has the same probability as $B$ started from $0$ hitting level $z-\mu$ by time  $x$. This can be calculated by the reflection principle; it equals $\Pr(|B(x)| > z-\mu)$. Since $\Pr(|B(x)| > a) \leq 2e^{-a^2/2x}$, the lemma follows.
\end{proof}

\begin{prop} \label{prop:SSdecay}
    Let $a \in \R$ and $0 < x_1, x_2 \leq T$. Set $|\mu| = \max_{y \in [0,x_2]} |b(y)|$. There is a constant $C = C_{a, T, |\mu|, x_2, n}$
    such that for all $u \in \R$ and $v \geq a$,
    $$ |S_{B,n,-x_1} \cdot S^{hypo(b)}_{B,n,x_2}(u,v) | \leq C (|u|^{n-1}+1) (|v|^n + 1)e^{-\frac{v^2}{4T}}.$$
\end{prop}

\begin{proof}
    Decompose $S_{B,n,-x_1} \cdot S^{hypo(b)}_{B,n,x_2}(u,v) = (I) + (II)$, where
    \begin{align*}
    (I) &= \int_{-\infty}^0 dz\, S_{B,n,-x_1}(u,z) S^{hypo(b)}_{B,n,x_2}(z,v) \\
    (II) &= \int_{0}^{\infty} dz\, S_{B,n,-x_1}(u,z) S^{hypo(b)}_{B,n,x_2}(z,v)
    \end{align*}

    Consider $(I)$ first. Since $z \leq 0$, and $b(0)=0$, the hitting time $\tau$ in $S^{hypo(b)}_{B,n,x_2}(z,v)$ is zero.
    So $S^{hypo(b)}_{B,n,x_2}(z,v) = \bar{S}_{B,n,x_2}(z,v)$. Then using Lemma \ref{lem:Hermitedecay} we find
    $$ |(I)| \leq \int_{-\infty}^{-v} dz\, C_{n,T} (|u-v-z|^{n-1}+1) (|z|^n +1) e^{- z^2/2x_2}.$$
    If $v \geq 0$ then the contribution to the integral above comes from the region $z \approx -v$. If $v < 0$ then it comes from the region $z \approx 0$.
    Since $v \geq a$ by assumption, we find that
    \begin{align*}
        |(I)| &\leq C'_{n,T} (|v|^m + 1) (|u|^{n-1}+1)\; [ e^{-v^2/2x_2} \ind{v \geq 0} + \ind{v \in [a,0]}] \\
        & \leq C_{a, n,T} (|v|^m + 1) (|u|^{n-1}+1) e^{-v^2/2T}
    \end{align*}

    Now consider $(II)$. Due to Lemma \ref{lem:Hermitedecay} we have
    $$ |S^{hypo(b)}_{B,n,x_2}(z,v)| \leq  C_{n,T} \, \E{ (|v-B(\tau)|^n+1) e^{\frac{(x_2-\tau)}{2}\partial^2}(B(\tau),v)\ind{\tau \leq x_2} \mid B(0)=z}.$$
    If $\tau \leq x_2$ then $|B(\tau)| \leq |\mu|$ because $B(\tau) = b(\tau)$ by continuity of $b$. Therefore,
    $$|S^{hypo(b)}_{B,n,x_2}(z,v)| \leq  C_{n,T,|\mu|} (|v|^n+1)\, \E{ e^{\frac{(x_2-\tau)}{2}\partial^2}(B(\tau),v)\ind{\tau \leq x_2} \mid B(0)=z}.$$
    Due to the strong Markov property of Brownian motion, we recognize the expectation above as $\Pr(\tau \leq x_2, B(x_2) \in dv \mid B(0)=z)$.
    By the Cauchy-Schwarz inequality,
    $$\Pr(\tau \leq x_2, B(x_2) \in dv \mid B(0)=z) \leq \Pr(\tau \leq x_2 \mid B(0)=z)^{1/2} \cdot \Pr(B(x_2) \in dv \mid B(0)=z)^{1/2}.$$
    By Lemma \ref{lem:tauprob} we have $\Pr(\tau \leq x_2 \mid B(0)=z) \leq \ind{z \leq \mu} + 2e^{-(z-\mu)^2/2x_2} \ind{z > \mu}$, where
    $\mu = \max_{y \in [0,x_2]} b(y)$. We also have $\Pr(B(x_2) \in dv \mid B(0)=z) = (\sqrt{2\pi x_2})^{-1/2} e^{-(z-v)^2/2x_2}$.
    Therefore,
    $$|S^{hypo(b)}_{B,n,x_2}(z,v)| \leq C_{n,T,|\mu|, x_2}(|v|^n+1) \left [ \ind{0\leq z \leq \mu} + e^{-\frac{(z-\mu)^2}{4x_2}} \ind{z > \mu}\right]
    e^{-\frac{(z-v)^2}{4x_2}}$$
    This estimate together with the estimate $|S_{B,n,-x_1}(u,z)| \leq C_{n,T} ( |u-z|^{n-1} + 1)$ implies
    $$ |(II)| \leq C_{n,T,|\mu|,x_2}\, (|u|^{n-1}+1)(|v|^n+1) e^{-v^2/4T}$$
    Combining the bounds on $(I)$ and $(II)$ proves the lemma.
\end{proof}

Proposition \ref{prop:SSdecay} shows that if we conjugate $S_{B,n,-x_i} \cdot S^{hypo(b)}_{B,n,x_j}$ by the factor in \eqref{eqn:conjugation}, then it is bounded in absolute value by $f_i(u) g_j(v)$ where $f_i$ is bounded over $\R$
and $g_j$ is integrable over $[a_j,\infty)$. So the Fredholm series expansion of $\dt{I-\chi_aK_B\chi_a}$ with the conjugation is absolutely convergent. This completes the proof of Theorem \ref{thm:Blpp}.

\subsection{Brownian tasep} \label{sec:btasep}
Consider the Brownian tasep \eqref{eqn:bbduality} with initial condition $X_n(0) = 0$ for every $n \geq 0$. We find that
$$ X_n(t) = - \max_{0 \leq k \leq n} \{ L_B(0,k;t,n)\} = -L_B(0,0;t,n) = -G(t,n)$$
In other words, $X$ is the space-time dual of $G$. Now looking at \eqref{eqn:GB}, if $b$ is the boundary condition for Brownian lpp, then the above shows that $X_0 = -b$ moves deterministically. Then $X_1$ is a Brownian motion reflected (to the left) off $-b$, $X_2$ is a Brownian motion reflected off $X_1$, etc.
Theorem \ref{thm:Blpp} thus implies
\begin{cor} \label{cor:Blpp}
   Consider Brownian tasep started from $X_n(0) = 0$ for $n \geq 0$.
    Let $X_n(t)$ be the location of particle number $n$ at time $t$. Let $X_0$ move deterministically according to a continuous function $-b$. The other particles move randomly as Brownian motion being reflected to the left off the particle to its right.
    Then for $n \geq 1$, distinct times $t_1, t_2, \ldots, t_k > 0$, and $a_1, \ldots, a_k \in \R$,
    \begin{equation} \label{eqn:btasepformula2}
    \Pr(X_n(t_i) \geq a_i; 1 \leq i \leq k) = \dt{I - \chi_{-a} J_B \chi_{-a}}_{L^2(\{1,2,\ldots,k\}\times \R)}
    \end{equation}
    where $\chi_{-a}(i,z) = \ind{z \geq -a_i}$ and the kernel $J_B$ is the kernel $K_B$ from Theorem \ref{thm:Blpp} with choice of parameters $x_i = t_i$ and $n$ the same. 
\end{cor}

\subsection{Brownian tasep with general initial conditions and Warren's formula} \label{sec:warrenformula}
Consider the Brownian tasep \eqref{eqn:bbduality} with general initial conditions $X_1(0) \geq X_2(0) \geq \cdots$. Warren \cite{Warren} provided a formula for the transition probability of this system. Let $x,y \in \R^N$ with $x_1 \geq x_2 \geq \cdots \geq x_N$ and likewise for $y$. Warren proved:
\begin{equation} \label{eqn:Warren}
    \Pr(X(t) \in dy \mid X(0) = x) = \dt{\partial^{j-i} \phi_t(y_j-x_i)}\,dy.
\end{equation}
Here $\phi_t$ is the probability density of a Normal random variable with mean 0 and variance $t$:
$$ \phi_t(x) = \frac{1}{\sqrt{2 \pi t}} e^{- \frac{x^2}{2t}}.$$
Actually, in Warren's convention, the particles are ordered in reverse. This corresponds to the negation of the process \eqref{eqn:bbduality}.
The formula does not change because $\phi_t$ is symmetric. Warren's formula is in fact the precursor to Johansson's formula \eqref{eqn:transition}.
It came about from the construction of the Warren process \cite{Warren}, which is a randomly evolving Gelfand-Tsetlin pattern in continuous space-time. Its rows perform Dyson Brownian motion subject to interlacing constraints and its right edge is Brownian tasep. The formula \eqref{eqn:Warren} comes from an intertwining relation between the bottom row and the right edge.

\subsubsection{Deriving Warren's formula}
Johansson's formula can be used to derive Warren's, which brings our discussion full circle.
In Brownian lpp, space and noise occupy continuous domains. As such, pushing and blocking mechanisms come together to become reflection.
One can thus derive \eqref{eqn:Warren} from \eqref{eqn:transition} by way of \eqref{eqn:Etransition}.

We will use Warren's convention for Brownian tasep, namely that particles are ordered from smallest to largest: $X_1(t) \leq X_2(t), \ldots, X_N(t)$.
We assume there are $N$ particles started from $x_1 = X_1(0) \leq \cdots \leq X_N(0) = x_N$.
Let $\ell > 0$ be a large scaling parameter. Consider the Exponential lpp model \eqref{eqn:Enoise} with boundary conditions
$G(m,0) = \sqrt{\ell}x_m + m$, $1 \leq m \leq N$.
Define, for $t \geq 0$,
$$ \vec{G}_{\ell}(t) = \left ( \frac{G(k,\lceil \ell t \rceil)-(\ell t+k)}{\sqrt{\ell}}; 1 \leq k \leq N\right)$$
Define for integers $x,k \geq 1$
$$ \eta_{\ell}^*(x,k) = \sum_{i=1}^x \omega_{k,i} \quad \eta_{\ell}^*(0,k)=0$$
Extend $\eta^*_{\ell}$ to $x \in [0,\infty)$ by linear interpolation.
We have that for every $k$,
$$ \bar{\eta}_{\ell}(t,k) = \frac{\eta^*_{\ell}(\ell t, k) - \ell t}{\sqrt{\ell}} \to B^*_k(t)$$
as $\ell \to \infty$, where $B^*_k$ are independent, standard Brownian motions.

We can write the corner growth function as
$$ G(k,n) = \max_{1 \leq m \leq k} \sqrt{\ell}x_m + m + L^*(1,m;k,n),$$
where
$$ L^*(1,m;k,n) = \max_{1 = x_{m-1} \leq x_m \leq x_{m+1} \leq \cdots \leq x_N = n} \sum_{i=m}^N \eta^*_{\ell}(x_i,i) - \eta^*_{\ell}(x_{i-1}-1,i).$$
Observe that
$$ \frac{L^*(1,m;k, \ell t) - (\ell t +k-m)}{\sqrt{\ell}} = \max_{0 = t_{m-1} \leq t_m \leq < \cdots \leq t_N=t}
\sum_{i=m}^N \bar{\eta}_{\ell}(t_i,i) - \bar{\eta}_{\ell}(t_{i-1},i) + O_p\big(\frac{N}{\sqrt{\ell}}\big).$$
As $\ell \to \infty$, the above converges in law to $L_B(0,m;t,N)$, where $L_B$ is the lpp function for the Brownian lpp model.
Hence,
$$\frac{G(k,\lceil \ell t \rceil) - (\ell t+k)}{\sqrt{\ell}} \to \max_{1 \leq m \leq k} x_m + L_B(0,m;t,N)$$
in law, for every $k$, as a process in $t$. Thus,
$$\vec{G}_{\ell}(t) \to X(t) = (X_1(t), \ldots, X_N(t)),$$
the latter being the Brownian tasep started from $x_1 \leq \cdots \leq x_N$.

In order to derive the transition probability \eqref{eqn:Warren}, we must compute
\begin{align*}
    \Pr(X(t) \in dy \mid X(0) = x)/dy  &= \lim_{\ell} \Pr(\vec{G}_{\ell}(t) \in dy \mid \vec{G}_{\ell}(0)=x)/dy \\
    &= \lim_{\ell} \Pr(G(k,\lceil \ell t \rceil) \in \ell t + k d(\sqrt{\ell y}))/ d(\sqrt{\ell}y) \; \ell^{N/2} \\
    &= \lim_{\ell} \dt{[\partial^{j-i}\cdot W_{\ell t}](\ell t + j + \sqrt{\ell}(y_j-x_i))} \ell^{N/2}.
\end{align*}
Observe that
$$[\partial^{j-i}\cdot W_{\ell t}](\ell t + j +\sqrt{\ell}(y_j-x_i)) = \partial_{y_j}^{j-i} [ W_{\ell t}(\ell t + j + \sqrt{\ell}(y_j-x_i)] \times \ell^{(i-j)/2}.$$
The last factor of $\ell^{(i-j)/2}$ is a conjugation and can be removed from the determinant.
Stirling's approximation of factorials : $n! = \sqrt{2 \pi n} (n/e)^n \times (1 + O(1/n))$ shows
$$W_{\ell t}(\ell t+j +\sqrt{\ell} (y_j-x_i)) = \frac{1}{\sqrt{2 \pi \ell t}} e^{- \frac{(y_j-x_i)^2}{2t}} \times (1 + O(1/\ell)).$$
Therefore, 
$$ \dt{[\partial^{j-i}\cdot W_{\ell t}](\ell t + j \sqrt{\ell}(y_j-x_i))} = \dt{[\partial^{j-i}\cdot\phi_t](y_j-x_i)} \ell^{-N/2} \times (1+O(1/\ell)).$$
So, finally,
$$\Pr(X(t) \in dy \mid X(0) = x)/dy = \dt{\partial^{j-i} \phi_t(y_j-x_i)}.$$
We also note that for $k \in \Z$,
$$ \partial^k \phi_t(x) = \frac{1}{2 \pi \mathbf{i}} \oint_{\Re(z) = \delta > 0} dz\, e^{\frac{t}{2}z^2 + xz} z^k.$$

\subsubsection{General initial conditions} 
One can use the formula \eqref{eqn:Warren} together with the Eynard-Mehta method and the solution of the corresponding bi-orthogonalization problem to express
$$ \Pr(X_{n_i}(t) \leq a_i, 1 \leq i \leq k \mid X(0)) = \dt{I-K}_{L^2(\{1,\ldots, k\} \times \R)}$$
as a Fredholm determinant of a kernel $K$. Another way to go about it is to use Theorem \ref{thm:Efixedtime}. We have that
$$\Pr(X_{n_i}(t) \leq a_i, 1 \leq i \leq k \mid X(0)) = \lim_{\ell \to \infty} \Pr\left(\frac{G(n_i, \lceil \ell t \rceil)-\ell t}{\sqrt{\ell}} \leq a_i; 1 \leq i \leq k\right)$$
Theorem \ref{thm:Efixedtime} provides a Fredholm determinant formula for the latter probability. On taking the large $\ell$ limit, one obtains a Fredholm determinant formula for Brownian tasep. This determinant formula has actually been derived in \cite{NQR} by taking a certain low-density limit of continuous time tasep. As such, we do not proceed further and end here.

\end{document}